\documentclass[12pt]{amsart}
\usepackage{amssymb,mathrsfs}

\begin{document}


\title[Invariant differential operators
on Siegel-Jacobi space]{Invariant differential operators
on \\ Siegel-Jacobi space }

\author{Jae-Hyun Yang}
\address{Department of Mathematics, Inha University, Incheon 402-751, Korea}
\email{jhyang@inha.ac.kr }


\newtheorem{theorem}{Theorem}[section]
\newtheorem{lemma}{Lemma}[section]
\newtheorem{proposition}{Proposition}[section]
\newtheorem{remark}{Remark}[section]
\newtheorem{definition}{Definition}[section]
\newtheorem{corollary}{Corollary}[section]

\renewcommand{\theequation}{\thesection.\arabic{equation}}
\renewcommand{\thetheorem}{\thesection.\arabic{theorem}}
\renewcommand{\thelemma}{\thesection.\arabic{lemma}}
\newcommand{\bbr}{\mathbb R}
\newcommand{\bbs}{\mathbb S}
\newcommand{\bn}{\bf n}
\newcommand\charf {\mbox{{\text 1}\kern-.24em {\text l}}}

\newcommand\BC{\mathbb C}
\newcommand\BZ{\mathbb Z}
\newcommand\BR{\Bbb R}
\newcommand\BQ{\mathbb Q}
\newcommand\Rmn{{\mathbb R}^{(m,n)}}

\newcommand\Rnn{{\mathbb R}^{(n,n)}}
\newcommand\Yd{{{\partial}\over {\partial Y}}}
\newcommand\Vd{{{\partial}\over {\partial V}}}
\newcommand\tr{\triangleright}
\newcommand\al{\alpha}
\newcommand\be{\beta}
\newcommand\g{\gamma}
\newcommand\gh{\Cal G^J}
\newcommand\G{\Gamma}
\newcommand\de{\delta}
\newcommand\e{\epsilon}
\newcommand\z{\zeta}
\newcommand\vth{\vartheta}
\newcommand\vp{\varphi}
\newcommand\om{\omega}
\newcommand\p{\pi}
\newcommand\la{\lambda}
\newcommand\lb{\lbrace}
\newcommand\lk{\lbrack}
\newcommand\rb{\rbrace}
\newcommand\rk{\rbrack}
\newcommand\s{\sigma}
\newcommand\w{\wedge}
\newcommand\fgj{{\frak g}^J}
\newcommand\lrt{\longrightarrow}
\newcommand\lmt{\longmapsto}
\newcommand\lmk{(\lambda,\mu,\kappa)}
\newcommand\Om{\Omega}
\newcommand\ka{\kappa}
\newcommand\ba{\backslash}
\newcommand\ph{\phi}
\newcommand\M{{\Cal M}}
\newcommand\bA{\bold A}
\newcommand\bH{\bold H}

\newcommand\Hom{\text{Hom}}
\newcommand\cP{\Cal P}
\newcommand\cH{\Cal H}

\newcommand\pa{\partial}

\newcommand\pis{\pi i \sigma}
\newcommand\sd{\,\,{\vartriangleright}\kern -1.0ex{<}\,}
\newcommand\wt{\widetilde}
\newcommand\fg{\frak g}
\newcommand\fk{\frak k}
\newcommand\fp{\frak p}
\newcommand\fs{\frak s}
\newcommand\fh{\frak h}
\newcommand\Cal{\mathcal}

\newcommand\fn{{\frak n}}
\newcommand\fa{{\frak a}}
\newcommand\fm{{\frak m}}
\newcommand\fq{{\frak q}}
\newcommand\CP{{\mathcal P}_g}
\newcommand\Hgh{{\mathbb H}_g \times {\mathbb C}^{(h,g)}}
\newcommand\BD{\mathbb D}
\newcommand\BH{\mathbb H}
\newcommand\CCF{{\mathcal F}_g}
\newcommand\CM{{\mathcal M}}
\newcommand\Ggh{\Gamma_{g,h}}
\newcommand\Chg{{\mathbb C}^{(h,g)}}

\newcommand\Ys{Y^{\ast}}
\newcommand\Vs{V^{\ast}}
\newcommand\LO{L_{\Omega}}
\newcommand\fac{{\frak a}_{\mathbb C}^{\ast}}
\newcommand\tra{ \textrm{tr}}

\thanks{\endgraf Keywords and phrases\,: invariants, invariant differential
operators, Siegel-Jacobi space.\\
\indent This work was supported by Basic Science Program through the National Research Foundation\\
\indent of Korea(NRF) funded by the Ministry of Education, Science and Technology (41493-01) and \\
\indent partially supported by the Max-Planck-Institut f{\"u}r Mathematik in Bonn.}

\begin{abstract} For two positive integers $m$ and $n$, we let ${\mathbb
H}_n$ be the Siegel upper half plane of degree $n$ and let
$\BC^{(m,n)}$ be the set of all $m\times n$ complex matrices. In
this article, we study differential operators on the
Siegel-Jacobi space ${\mathbb H}_n\times {\mathbb C}^{(m,n)}$ that
are invariant under the {\it natural} action of the Jacobi group
$Sp(n,\BR)\ltimes H_{\BR}^{(n,m)}$ on ${\mathbb H}_n\times
{\mathbb C}^{(m,n)}$, where $H_{\BR}^{(n,m)}$ denotes the
Heisenberg group. We give some explicit invariant differential operators.
We present important problems which are natural. We give some partial solutions for these natural problems.
\end{abstract}
\maketitle

\begin{section}{{\large\bf Introduction}}
\setcounter{equation}{0}

For a given fixed positive integer $n$, we let
$${\mathbb H}_n=\,\{\,\Omega\in \BC^{(n,n)}\,|\ \Om=\,^t\Om,\ \ \ \text{Im}\,\Om>0\,\}$$
be the Siegel upper half plane of degree $n$ and let
$$Sp(n,\BR)=\{ M\in \BR^{(2n,2n)}\ \vert \ ^t\!MJ_nM= J_n\ \}$$
be the symplectic group of degree $n$, where $F^{(k,l)}$ denotes
the set of all $k\times l$ matrices with entries in a commutative
ring $F$ for two positive integers $k$ and $l$, $^t\!M$ denotes
the transpose matrix of a matrix $M$ and
$$J_n=\begin{pmatrix} 0&I_n\\
                   -I_n&0\end{pmatrix}.$$
$Sp(n,\BR)$ acts on $\BH_n$ transitively by
\begin{equation}
M\cdot\Om=(A\Om+B)(C\Om+D)^{-1},
\end{equation} where $M=\begin{pmatrix} A&B\\
C&D\end{pmatrix}\in Sp(n,\BR)$ and $\Om\in \BH_n.$

 \vskip 0.1cm For two
positive integers $m$ and $n$, we consider the Heisenberg group
$$H_{\BR}^{(n,m)}=\big\{\,(\la,\mu;\ka)\,|\ \la,\mu\in \BR^{(m,n)},\ \kappa\in\BR^{(m,m)},\
\ka+\mu\,^t\la\ \text{symmetric}\ \big\}$$ endowed with the
following multiplication law
$$\big(\la,\mu;\ka\big)\circ \big(\la',\mu';\ka'\big)=\big(\la+\la',\mu+\mu';\ka+\ka'+\la\,^t\mu'-
\mu\,^t\la'\big)$$
with $\big(\la,\mu;\ka\big),\big(\la',\mu';\ka'\big)\in H_{\BR}^{(n,m)}.$
We define the semidirect product of
$Sp(n,\BR)$ and $H_{\BR}^{(n,m)}$
$$G^J=Sp(n,\BR)\ltimes H_{\BR}^{(n,m)}$$
endowed with the following multiplication law
$$
\big(M,(\lambda,\mu;\kappa)\big)\cdot\big(M',(\lambda',\mu';\kappa'\,)\big)
=\, \big(MM',(\tilde{\lambda}+\lambda',\tilde{\mu}+ \mu';
\kappa+\kappa'+\tilde{\lambda}\,^t\!\mu'
-\tilde{\mu}\,^t\!\lambda'\,)\big)$$ with $M,M'\in Sp(n,\BR),
(\lambda,\mu;\kappa),\,(\lambda',\mu';\kappa') \in
H_{\BR}^{(n,m)}$ and
$(\tilde{\lambda},\tilde{\mu})=(\lambda,\mu)M'$. Then $G^J$ acts
on $\BH_n\times \BC^{(m,n)}$ transitively by
\begin{equation}
\big(M,(\lambda,\mu;\kappa)\big)\cdot
(\Om,Z)=\Big(M\cdot\Om,(Z+\lambda \Om+\mu)
(C\Omega+D)^{-1}\Big), \end{equation} where $M=\begin{pmatrix} A&B\\
C&D\end{pmatrix} \in Sp(n,\BR),\ (\lambda,\mu; \kappa)\in
H_{\BR}^{(n,m)}$ and $(\Om,Z)\in \BH_n\times \BC^{(m,n)}.$ We note
that the Jacobi group $G^J$ is {\it not} a reductive Lie group and
that the homogeneous space ${\mathbb H}_n\times \BC^{(m,n)}$ is not a
symmetric space. We refer to \cite{BS,EZ,YJH1,YJH2,YJH3,YJH4,YJH5,YJH6,YJH7,YJH8,YJH9} about
automorphic forms on $G^J$ and topics related to the content of
this paper. From now on, for brevity we write
$\BH_{n,m}=\BH_n\times \BC^{(m,n)},$ called the Siegel-Jacobi space of degree $n$ and index $m$.

\vskip 0.2cm The aim of this paper is to study differential
operators on ${\mathbb H}_{n,m}$ which are invariant under the
{\it natural} action (1.2) of $G^J$. The study of these invariant differential operators on the
Siegel-Jacobi space $\BH_{n,m}$ is interesting and important in the aspects of invariant theory,
arithmetic and geometry.
This article is organized as follows. In
Section 2, we review differential operators on ${\mathbb H}_n$
invariant under the action (1.1) of $Sp(n,\BR)$. We let $\BD ({\mathbb H}_n)$
denote the algebra of all differential operators on ${\mathbb H}_n$ that are invariant under the action
(1.1). According to the work of Harish-Chandra \cite{HC1,HC2}, we see that $\BD ({\mathbb H}_n)$ is a commutative algebra which is
isomorphic to the center of the universal enveloping algebra of the complexification of the Lie algebra of $Sp(n,\BR)$.
We briefly describe the work of Maass\,\cite{M2} about constructing explicit algebraically independent generators
of $\BD ({\mathbb H}_n)$ and Shimura's construction\,\cite{Sh} of canonically defined algebraically independent generators
of $\BD ({\mathbb H}_n)$.
In Section 3, we
study differential operators on ${\mathbb H}_{n,m}$
invariant under the action (1.2) of $G^J$. For two positive integers $m$ and $n$, we let
\begin{equation*}
T_{n,m}=\,\left\{\, (\omega,z)\,|\ \omega=\,{}^t\omega\in \BC^{(n,n)},\ z\in\BC^{(m,n)}\,\right\}
\end{equation*}
be the complex vector space of dimension ${{n(n+1)}\over 2}+mn.$
From the adjoint action of the Jacobi group $G^J$,
we have the {\it natural action} of the unitary group $U(n)$
on $T_{n,m}$ given by
\begin{equation}
u\cdot (\omega,z)=\,( u\,\omega\,{}^t u, z\,{}^tu),\quad u\in U(n),\ (\omega,z)\in T_{n,m}.
\end{equation}
The action (1.3) of $U(n)$ induces canonically the representation $\tau$ of $U(n)$ on the polynomial algebra
$\textrm{Pol}(T_{n,m})$ consisting of complex valued polynomial functions on $T_{n,m}.$
Let $ \textrm{Pol}(T_{n,m})^{U(n)}$ denote the subalgebra of $ \textrm{Pol}(T_{n,m})$ consisting of all polynomials on
$T_{n,m}$ invariant under the representation $\tau$ of $U(n)$, and
$\BD ({\mathbb H}_{n,m})$ denote the algebra of all differential operators on ${\mathbb H}_{n,m}$ invariant under the
action $(1.2)$ of $G^J$. We see that there is a canonically defined linear bijection of $\textrm{Pol}(T_{n,m})^{U(n)}$
onto $\BD ({\mathbb H}_{n,m})$ which is not multiplicative. We will see that $\BD ({\mathbb H}_{n,m})$ is {\it not} commutative.
The main important problem is to find explicit generators of $ \textrm{Pol}(T_{n,m})^{U(n)}$
and explicit generators of $\BD ({\mathbb H}_{n,m})$. We propose several natural problems.
We want to mention that at this moment it is
quite complicated and difficult to find the explicit generators of $\BD ({\mathbb H}_{n,m})$ and to express invariant differential operators on
$\BH_{n,m}$ explicitly.
In Section 4, we gives some examples of explicit $G^J$-invariant differential operators on ${\mathbb H}_{n,m}$
that are obtained by complicated calculations.
In Section 5, we deal with the special case $n=m=1$ in detail. We give complete solutions of the problems that are
proposed in Section 3.
In Section 6, we deal with the case that $n=1$ and $m$ is arbitrary. We give some partial solutions for the problems proposed
in Section 3.
In the final section, using these invariant differential operators on the Siegel-Jacobi space,
we discuss a notion of Maass-Jacobi forms.

\vskip 0.31cm \noindent {\bf Acknowledgements:} This work was in part done during my stay at the Max-Planck-Institut f{\"u}r Mathematik in
Bonn. I am very grateful for the hospitality and financial support. I also thank the National Research Foundation of Korea for its financial support.
Finally I would like to give my hearty thanks to Don Zagier, Eberhard Freitag, Rainer Weissauer, Hiroyuki Ochiai and Minoru Itoh for their
interests in this work and fruitful discussions.

\vskip 0.31cm \noindent {\bf Notations:} \ \ We denote by
$\BQ,\,\BR$ and $\BC$ the field of rational numbers, the field of
real numbers and the field of complex numbers respectively. We
denote by $\BZ$ and $\BZ^+$ the ring of integers and the set of
all positive integers respectively. The symbol ``:='' means that
the expression on the right is the definition of that on the left.
For two positive integers $k$ and $l$, $F^{(k,l)}$ denotes the set
of all $k\times l$ matrices with entries in a commutative ring
$F$. For a square matrix $A\in F^{(k,k)}$ of degree $k$,
$\textrm{tr}(A)$ denotes the trace of $A$. For any $M\in F^{(k,l)},\
^t\!M$ denotes the transpose matrix of $M$. $I_n$ denotes the
identity matrix of degree $n$. For $A\in F^{(k,l)}$ and $B\in
F^{(k,k)}$, we set $B[A]=\,^tABA.$ For a complex matrix $A$,
${\overline A}$ denotes the complex {\it conjugate} of $A$. For
$A\in \BC^{(k,l)}$ and $B\in \BC^{(k,k)}$, we use the abbreviation
$B\{ A\}=\,^t{\overline A}BA.$ For a positive integer $n$, $I_n$ denotes the
identity matrix of degree $n$.

\end{section}

\vskip 1cm

\begin{section}{{\large\bf Invariant Differential Operators on the Siegel Space }}
\setcounter{equation}{0}

\newcommand\POB{ {{\partial\ }\over {\partial{\overline \Omega}}} }
\newcommand\PZB{ {{\partial\ }\over {\partial{\overline Z}}} }
\newcommand\PX{ {{\partial\ }\over{\partial X}} }
\newcommand\PY{ {{\partial\ }\over {\partial Y}} }
\newcommand\PU{ {{\partial\ }\over{\partial U}} }
\newcommand\PV{ {{\partial\ }\over{\partial V}} }
\newcommand\PO{ {{\partial\ }\over{\partial \Omega}} }
\newcommand\PZ{ {{\partial\ }\over{\partial Z}} }

\newcommand\POBS{ {{\partial\ \,}\over {\partial{\overline \Omega}_*} } }
\newcommand\PZBS{ {{\partial\ \,}\over {\partial{\overline Z_*}}} }
\newcommand\PXS{ {{\partial\ \,}\over{\partial X_*}} }
\newcommand\PYS{ {{\partial\ \,}\over {\partial Y_*}} }
\newcommand\PUS{ {{\partial\ \,}\over{\partial U_*}} }
\newcommand\PVS{ {{\partial\ \,}\over{\partial V_*}} }
\newcommand\POS{ {{\partial\ \,}\over{\partial \Omega_*}} }
\newcommand\PZS{ {{\partial\ \,}\over{\partial Z_*}} }

\vskip 0.21cm For a coordinate $\Om=(\omega_{ij})\in\BH_n,$ we write
$\Om=X+i\,Y$ with $X=(x_{ij}),\ Y=(y_{ij})$ real. We put
$d\Om=\big(d\om_{ij}\big)$ and $d{\overline\Om}=\big(d{\overline
{\om}}_{ij}\big)$. We also put
$$\PO=\,\left(\,
{ {1+\delta_{ij}}\over 2}\, { {\partial\ \,}\over {\partial \om_{ij} }
} \,\right) \qquad\text{and}\qquad \POB=\,\left(\, {
{1+\delta_{ij}}\over 2}\, { {\partial\ \,}\over {\partial {\overline
{\om}}_{ij} } } \,\right).$$ Then for a positive real number $A$,
\begin{equation}
ds_{n;A}^2=A\, \textrm{tr} \Big(Y^{-1}d\Om\,
Y^{-1}d{\overline\Om}\,\Big)\end{equation} is a
$Sp(n,\BR)$-invariant K{\"a}hler metric on $\BH_n$
(cf.\,\cite{Si1, Si2}), where $ \textrm{tr}(M)$ denotes the trace of a
square matrix $M$. H. Maass \cite{M} proved that the Laplacian of
$ds^2_{n;A}$ is given by
\begin{equation}
\Delta_{n;A}=\,{\frac 4A}\, \textrm{tr} \left( Y\,\,
{}^{{}^{{}^{{}^\text{\scriptsize $t$}}}}\!\!\!
\left(Y\POB\right)\PO\right).\end{equation} And
\begin{equation*}
dv_n(\Om)=(\det Y)^{-(n+1)}\prod_{1\leq i\leq j\leq n}dx_{ij}\,
\prod_{1\leq i\leq j\leq n}dy_{ij}\end{equation*} is a
$Sp(n,\BR)$-invariant volume element on
$\BH_n$\,(cf.\,\cite[p.\,130]{Si2}).

\vskip 0.2cm For brevity, we write $G=Sp(n,\BR).$ The isotropy
subgroup $K$ at $iI_n$ for the action (1.1) is a maximal compact
subgroup given by
\begin{equation*}
K=\left\{ \begin{pmatrix} A & -B \\ B & A \end{pmatrix} \Big| \
A\,^t\!A+ B\,^t\!B=I_n,\ A\,^t\!B=B\,^t\!A,\ A,B\in
\BR^{(n,n)}\,\right\}.
\end{equation*}

\noindent Let $\fk$ be the Lie algebra of $K$. Then the Lie
algebra $\fg$ of $G$ has a Cartan decomposition $\fg=\fk\oplus
\fp$, where
\begin{equation*}
\frak g=\left\{ \begin{pmatrix} X_1 & \ \ X_2 \\ X_3 & -\,{}^tX_1
\end{pmatrix}\,\Big|\ X_1,X_2,X_3\in\BR^{(n,n)},\ X_2=\,{}^tX_2,\ X_3=\,{}^tX_3\,\right\},
\end{equation*}

\begin{equation*}
\frak k=\left\{ \begin{pmatrix} X &  -Y \\ Y & \ X
\end{pmatrix}\in\BR^{(2n,2n)}\,\Big|
\ {}^tX+X=0,\ Y=\,{}^tY\,\right\},
\end{equation*}

\begin{equation*}
\fp=\left\{ \begin{pmatrix} X & \ Y \\ Y & -X \end{pmatrix} \Big| \
X=\,^tX,\ Y=\,^tY,\ X,Y\in \BR^{(n,n)}\,\right\}.
\end{equation*}

The subspace $\fp$ of $\fg$ may be regarded as the tangent space
of $\BH_n$ at $iI_n.$ The adjoint representation of $G$ on $\fg$
induces the action of $K$ on $\fp$ given by
\begin{equation}
k\cdot Z=\,kZ\,^tk,\quad k\in K,\ Z\in \fp.
\end{equation}

Let $T_n$ be the vector space of $n\times n$ symmetric complex
matrices. We let $\Psi: \fp\lrt T_n$ be the map defined by
\begin{equation}
\Psi\left( \begin{pmatrix} X & \ Y \\ Y & -X \end{pmatrix}
\right)=\,X\,+\,i\,Y, \quad \begin{pmatrix} X & \ Y \\ Y & -X
\end{pmatrix}\in \fp.
\end{equation}

\noindent We let $\delta:K\lrt U(n)$ be the isomorphism defined by
\begin{equation}
\delta\left( \begin{pmatrix} A & -B \\ B & A \end{pmatrix}
\right)=\,A\,+\,i\,B, \quad \begin{pmatrix} A & -B \\ B & A
\end{pmatrix}\in K,
\end{equation}

\noindent where $U(n)$ denotes the unitary group of degree $n$. We
identify $\fp$ (resp. $K$) with $T_n$ (resp. $U(n)$) through the
map $\Psi$ (resp. $\delta$). We consider the action of $U(n)$ on
$T_n$ defined by
\begin{equation}
h\cdot \omega=\,h\omega \,^th,\quad h\in U(n),\ \omega\in T_n.
\end{equation}

\noindent Then the adjoint action (2.3) of $K$ on $\fp$ is
compatible with the action (2.6) of $U(n)$ on $T_n$ through the
map $\Psi.$ Precisely for any $k\in K$ and $Z\in \fp$, we get
\begin{equation}
\Psi(k\,Z \,^tk)=\delta(k)\,\Psi(Z)\,^t\delta (k).
\end{equation}

\noindent The action (2.6) induces the action of $U(n)$ on the
polynomial algebra $ \textrm{Pol}(T_n)$ and the symmetric algebra
$S(T_n)$ respectively. We denote by $ \textrm{Pol}(T_n)^{U(n)}$
$\Big( \textrm{resp.}\ S(T_n)^{U(n)}\,\Big)$ the subalgebra of $
\textrm{Pol}(T_n)$ $\Big( \textrm{resp.}\ S(T_n)\,\Big)$
consisting of $U(n)$-invariants. The following inner product $(\
,\ )$ on $T_n$ defined by $$(Z,W)= \, \textrm{tr}
\big(Z\,{\overline W}\,\big),\quad Z,W\in T_n$$

\noindent gives an isomorphism as vector spaces
\begin{equation}
T_n\cong T_n^*,\quad Z\mapsto f_Z,\quad Z\in T_n,
\end{equation}

\noindent where $T_n^*$ denotes the dual space of $T_n$ and $f_Z$
is the linear functional on $T_n$ defined by
$$f_Z(W)=(W,Z),\quad W\in T_n.$$

\noindent It is known that there is a canonical linear bijection
of $S(T_n)^{U(n)}$ onto the algebra ${\mathbb D}(\BH_n)$ of
differential operators on $\BH_n$ invariant under the action (1.1)
of $G$. Identifying $T_n$ with $T_n^*$ by the above isomorphism
(2.8), we get a canonical linear bijection
\begin{equation}
\Theta_n:\textrm{Pol}(T_n)^{U(n)} \lrt {\mathbb D}(\BH_n)
\end{equation}

\noindent of $ \textrm{Pol}(T_n)^{U(n)}$ onto ${\mathbb
D}(\BH_n)$. The map $\Theta_n$ is described explicitly as follows.
Similarly the action (2.3) induces the action of $K$ on the
polynomial algebra $ \textrm{Pol}(\fp)$ and the symmetric algebra $S(\fp)$ respectively.
Through the map $\Psi$, the subalgebra $ \textrm{Pol}(\fp)^K$ of $
\textrm{Pol}(\fp)$ consisting of $K$-invariants is isomorphic to $
\textrm{Pol}(T_n)^{U(n)}$. We put $N=n(n+1)$. Let $\left\{
\xi_{\alpha}\,|\ 1\leq \alpha \leq N\, \right\}$ be a basis of a real vector space
$\fp$. If $P\in \textrm{Pol}(\fp)^K$, then
\begin{equation}
\Big(\Theta_n (P)f\Big)(gK)=\left[ P\left( {{\partial\ }\over {\partial
t_{\al}}}\right)f\left(g\,\text{exp}\, \left(\sum_{\al=1}^N
t_{\al}\xi_{\al}\right) K\right)\right]_{(t_{\al})=0},
\end{equation} where $f\in C^{\infty}({\mathbb H}_{n})$. We refer to \cite{He1,He2} for more detail. In
general, it is hard to express $\Phi(P)$ explicitly for a
polynomial $P\in \textrm{Pol}(\fp)^K$.

\vskip 0.3cm According to the work of Harish-Chandra \cite{HC1,HC2}, the
algebra ${\mathbb D}(\BH_n)$ is generated by $n$ algebraically
independent generators and is isomorphic to the commutative ring
$\BC [x_1,\cdots,x_n]$ with $n$ indeterminates. We note that $n$
is the real rank of $G$. Let $\fg_{\BC}$ be the complexification
of $\fg$. It is known that $\BD(\BH_n)$ is isomorphic to the
center of the universal enveloping algebra of $\fg_{\BC}$.

\vskip 0.3cm Using a classical invariant theory (cf.\,\cite{Ho, W},
we can show that $\textrm{Pol}(T_n)^{U(n)}$ is generated by the
following algebraically independent polynomials
\begin{equation}
q_j (\omega)=\,\textrm{tr}\Big( \big(\omega {\overline
\omega}\big)^j\,\Big),\quad \omega\in T_n, \quad j=1,2,\cdots,n.
\end{equation}

For each $j$ with $1\leq j\leq n,$ the image $\Theta_n(q_j)$ of $q_j$
is an invariant differential operator on $\BH_n$ of degree $2j$.
The algebra ${\mathbb D}(\BH_n)$ is generated by $n$ algebraically
independent generators $\Theta_n(q_1),\Theta_n(q_2),\cdots,\Theta_n(q_n).$ In
particular,
\begin{equation}
\Theta_n(q_1)=\,c_1\, \textrm{tr}\! \left( Y\,
{}^{{}^{{}^{{}^\text{\scriptsize $t$}}}}\!\!\!
\left(Y\POB\right)\!\PO\right)\quad  \textrm{for\ some
constant}\ c_1.
\end{equation}

\noindent We observe that if we take $\omega=x+\,i\,y\in T_n$ with real $x,y$,
then $q_1(\omega)=q_1(x,y)=\,\textrm{tr}\big( x^2 +y^2\big)$ and
\begin{equation*}
q_2(\omega)=q_2(x,y)=\, \textrm{tr}\Big(
\big(x^2+y^2\big)^2+\,2\,x\big(xy-yx)y\,\Big).
\end{equation*}

\vskip 0.3cm It is a natural question to express the images
$\Theta_n(q_j)$ explicitly for $j=2,3,\cdots,n.$ We hope that the images $\Theta_n(q_j)$ for
$j=2,3,\cdots,n$ are expressed in the form of the $\textit{trace}$
as $\Phi(q_1)$.

\vskip 0.3cm H. Maass \cite{M2} found algebraically independent generators $H_1,H_2,\cdots,H_n$ of ${\mathbb D}(\BH_n)$.
We will describe $H_1,H_2,\cdots,H_n$ explicitly. For $M=\begin{pmatrix} A&B\\
C&D\end{pmatrix} \in Sp(n,\BR)$ and $\Omega=X+iY\in \BH_n$ with real $X,Y$, we set
\begin{equation*}
\Omega_*=\,M\!\cdot\!\Omega=\,X_*+\,iY_*\quad \textrm{with}\ X_*,Y_*\ \textrm{real}.
\end{equation*}
We set
\begin{eqnarray*}
K&=&\,\big( \Omega-{\overline\Omega}\,\big)\PO=\,2\,i\,Y \PO,\\
\Lambda&=&\,\big( \Omega-{\overline\Omega}\,\big)\POB=\,2\,i\,Y \POB,\\
K_*&=& \,\big( \Omega_*-{\overline\Omega}_*\,\big)\POS=\,2\,i\,Y_* \POS,\\
\Lambda_*&=&\,\big( \Omega_*-{\overline\Omega}_*\,\big)\POBS=\,2\,i\,Y_* \POBS.
\end{eqnarray*}
Then it is easily seen that
\begin{equation}
K_*=\,{}^t(C{\overline\Om}+D)^{-1}\,{}^t\!\left\{ (C\Omega+D)\,{}^t\!K \right\},
\end{equation}

\begin{equation}
\Lambda_*=\,{}^t(C{\Om}+D)^{-1}\,{}^t\!\left\{ (C{\overline\Omega}+D)\,{}^t\!\Lambda \right\}
\end{equation}
and
\begin{equation}
{}^t\!\left\{ (C{\overline\Omega}+D)\,{}^t\!\Lambda \right\}=\,\Lambda\,{}^t(C{\overline\Omega}+D)
-{{n+1}\over 2} \,\big( \Omega-{\overline\Omega}\,\big)\,{}^t\!C.
\end{equation}

Using Formulas (2.13),\,(2.14) and (2.15), we can show that

\begin{equation}
\Lambda_*K_* \,+\,{{n+1}\over 2}K_*=\,{}^t(C{\Om}+D)^{-1}\,{}^{{}^{{}^{{}^\text{\scriptsize $t$}}}}\!\!\!
\left\{ (C{\Omega}+D)\,{}^{{}^{{}^{{}^\text{\scriptsize $t$}}}}\!\!\!
\left( \Lambda K \,+\, {{n+1}\over 2}K \right)\right\}.
\end{equation}

\noindent Therefore we get
\begin{equation}
\textrm{tr}\!  \left( \Lambda_*K_* \,+\,{{n+1}\over 2}K_* \right) =\,
 \textrm{tr}\!   \left( \Lambda K \,+\, {{n+1}\over 2}K \right).
\end{equation}
We set
\begin{equation}
A^{(1)}=\,\Lambda K \,+\, {{n+1}\over 2}K .
\end{equation}

We define $A^{(j)}\,(j=2,3,\cdots,n)$ recursively by
\begin{eqnarray}
A^{(j)}&=&\, A^{(1)}A^{(j-1)}- {{n+1}\over 2}\,\Lambda\, A^{(j-1)}\,+\,{\frac 12}\,\Lambda\,\textrm{tr}\!\left(
A^{(j-1)} \right)\\
& & \ \ \,+\, {\frac 12}\,\big( \Omega-{\overline\Omega}\,\big) \,
{}^{{}^{{}^\text{\scriptsize $t$}}}\!\!\!\left\{
\big( \Omega-{\overline\Omega}\,\big)^{-1}
\,{}^t\!\left( \,{}^t\!\Lambda\,{}^t\!A^{(j-1)}\right)\right\}.\nonumber
\end{eqnarray}

\noindent We set
\begin{equation}
H_j=\,\textrm{tr}\!  \left( A^{(j)} \right),\quad j=1,2,\cdots,n.
\end{equation}
As mentioned before, Maass proved that $H_1,H_2,\cdots,H_n$ are algebraically independent generators
of ${\mathbb D}(\BH_n)$.

\vskip 0.3cm In fact, we see that
\begin{equation}
 -H_1 =\Delta_{n;1}=\,4\, \textrm{tr}\! \left( Y\,
{}^{{}^{{}^{{}^\text{\scriptsize $t$}}}}\!\!\!
\left(Y\POB\right)\!\PO\right).
\end{equation}
is the Laplacian for the invariant metric $ds^2_{n;1}$ on $\BH_n$.

\vskip 0.3cm
\noindent
{\bf Conjecture.} For $j=2,3,\cdots,n,\ \Theta_n(q_j)=\,c_j\,H_j$ for a suitable constant $c_j.$

\vskip 0.53cm\noindent $ \textbf{Example 2.1.}$ We consider the
case $n=1.$ The algebra $ \textrm{Pol}(T_1)^{U(1)}$ is generated
by the polynomial
\begin{equation*}
q(\omega)=\omega\,{\overline \omega},\quad \omega=x+ \,iy\in \BC \ \textrm{with}\ x,y \ \textrm{real}.
\end{equation*}

Using Formula (2.10), we get

\begin{equation*}
\Theta_1 (q)=\,4\,y^2 \left( { {\partial^2}\over {\partial x^2} }+{
{\partial^2}\over {\partial y^2} }\,\right).
\end{equation*}

\noindent Therefore $\BD (\BH_1)=\BC\big[ \Theta_1(q)\big]=\,\BC [H_1].$

\vskip 0.3cm\noindent $ \textbf{Example 2.2.}$ We consider the
case $n=2.$ The algebra $ \textrm{Pol}(T_2)^{U(2)}$ is generated
by the polynomial
\begin{equation*}
q_1(\omega)=\,\tra \big(\omega\,{\overline \omega}\,\big),\quad q_2(\omega)=\,\tra
\Big( \big(\omega\,{\overline \omega}\big)^2\Big), \quad \omega\in T_2.
\end{equation*}

Using Formula (2.10), we may express $\Theta_2(q_1)$ and $\Theta_2(q_2)$
explicitly. $\Theta_2 (q_1)$ is expressed by Formula (2.12). The
computation of $\Theta_2(q_2)$ might be quite tedious. We leave the
detail to the reader. In this case, $\Theta_2 (q_2)$ was essentially
computed in \cite{BC}, Proposition 6. Therefore
\begin{equation*}
\BD (\BH_2)=\BC\big[
\Theta_2(q_1), \Theta_2(q_2)\big]=\,\BC [H_1,H_2].
\end{equation*}
In fact, the center of the universal enveloping algebra
${\mathscr U}(\fg_{\BC})$ was computed in \cite{BC}.

\vskip 0.3cm G. Shimura \cite{Sh} found canonically defined algebraically independent generators
of $\BD(\BH_n)$. We will describe his way of constructing those generators roughly. Let $K_\BC,\,
{\frak g}_\BC,\,{\mathfrak k}_\BC,{\mathfrak p}_\BC,\cdots$ denote the complexication of $K,\,{\mathfrak g},\,
{\mathfrak k},\,{\mathfrak p},\cdots$ respectively. Then we have the Cartan decomposition
\begin{equation*}
{\frak g}_\BC=\,{\mathfrak k}_\BC + {\mathfrak p}_\BC,\quad
{\mathfrak p}_\BC=\,{\mathfrak p}_\BC^+ + {\mathfrak p}_\BC^-
\end{equation*}
with the properties
\begin{equation*}
[{\mathfrak k}_\BC,{\mathfrak p}_\BC^{\pm}]\subset {\mathfrak p}_\BC^{\pm},\ \ \
[{\mathfrak p}_\BC^{+},{\mathfrak p}_\BC^+]=[{\mathfrak p}_\BC^-,{\mathfrak p}_\BC^-]=\{0\},
\ \ \ [{\mathfrak p}_\BC^+,{\mathfrak p}_\BC^-]=\,{\mathfrak k}_\BC,
\end{equation*}
where
\begin{equation*}
{\frak g}_\BC=\left\{ \begin{pmatrix} X_1 & \ \ X_2 \\ X_3 & -\,{}^tX_1
\end{pmatrix}\,\Big|\ X_1,X_2,X_3\in\BC^{(n,n)},\ X_2=\,{}^tX_2,\ X_3=\,{}^tX_3\,\right\},
\end{equation*}
\begin{equation*}
{\frak k}_\BC=\left\{ \begin{pmatrix} A &  -B \\ B & \ A
\end{pmatrix}\in\BC^{(2n,2n)}\,\Big|
\ {}^tA+A=0,\ B=\,{}^tB\,\right\},
\end{equation*}
\begin{eqnarray*}
{\mathfrak p}_\BC&=&\,\left\{ \begin{pmatrix} X & \ Y \\ Y & -X
\end{pmatrix}\in\BC^{(2n,2n)}\,\Big|
\ X=\,{}^tX,\ Y=\,{}^tY\,\right\},\\
{\mathfrak p}_\BC^+&=&\,\left\{ \begin{pmatrix} Z & iZ \\ iZ & -Z
\end{pmatrix}\in\BC^{(2n,2n)}\,\Big|
\ Z=\,{}^tZ\in \BC^{(n,n)}\,\right\},\\
{\mathfrak p}_\BC^-&=&\,\left\{ \begin{pmatrix} \ Z & -iZ \\ -iZ & \,-Z
\end{pmatrix}\in\BC^{(2n,2n)}\,\Big|
\ Z=\,{}^tZ\in \BC^{(n,n)}\,\right\}.
\end{eqnarray*}

For a complex vector space $W$ and a nonnegative integer $r$, we denote by $ \textrm{Pol}_r(W)$
the vector space of complex-valued homogeneous polynomial functions on $W$ of degree $r$. We put
$$\textrm{Pol}^r(W):=\sum_{s=0}^r \textrm{Pol}_s(W).$$
$ \textrm{Ml}_r(W)$ denotes the vector space of all $\BC$-multilinear maps of $W\times \cdots
\times W\,( r \ \textrm{copies})$ into $\BC$.
An element $Q$ of $ \textrm{Ml}_r(W)$ is called {\it symmetric} if
$$Q(x_1,\cdots,x_r)=\,Q(x_{\pi(1)},\cdots,x_{\pi(r)})$$ for each permutation $\pi$ of
$\{ 1,2,\cdots,r\}.$ Given $P\in \textrm{Pol}_r(W)$, there is a unique element symmetric element
$P_*$ of $ \textrm{Ml}_r(W)$ such that
\begin{equation}
P(x)=\,P_*(x,\cdots,x)\qquad \textrm{for all}\ x\in W.
\end{equation}
Moreover the map $P\mapsto P_*$ is a $\BC$-linear bijection of $\textrm{Pol}_r(W)$ onto the
set of all symmetric elements of $ \textrm{Ml}_r(W)$. We let $S_r(W)$ denote the subspace
consisting of all homogeneous elements of degree $r$ in the symmetric algebra $S(W)$. We note that
$ \textrm{Pol}_r(W)$ and $S_r(W)$ are dual to each other with respect to the pairing
\begin{equation}
\langle \alpha,x_1\cdots x_r \rangle=\,\alpha_*(x_1,\cdots,x_r)\qquad (x_i\in W,\ \alpha
\in \textrm{Pol}_r(W)).
\end{equation}

\vskip 0.2cm
Let ${\mathfrak p}_\BC^*$ be the dual space of ${\mathfrak p}_\BC$, that is, ${\mathfrak p}_\BC^*=
\textrm{Pol}_1({\mathfrak p}_\BC).$ Let $\{ X_1,\cdots,X_N\}$ be a basis of ${\mathfrak p}_\BC$ and
$\{ Y_1,\cdots,Y_N\}$ be the basis of ${\mathfrak p}_\BC^*$ dual to $\{ X_\nu\},$ where
$N=n(n+1)$. We note that $\textrm{Pol}_r({\mathfrak p}_\BC)$ and
$\textrm{Pol}_r({\mathfrak p}_\BC^*)$ are dual to each other with respect to the pairing
\begin{equation}
\langle\alpha,\beta\rangle=\sum \alpha_*(X_{i_1},\cdots,X_{i_r})\,\beta_* (Y_{i_1},\cdots,Y_{i_r}),\
\end{equation}
where $\alpha\in \textrm{Pol}_r({\mathfrak p}_\BC),\ \beta\in \textrm{Pol}_r({\mathfrak p}_\BC^*)$ and
$(i_1,\cdots,i_r)$ runs over $\{1,\cdots,N\}^r.$ Let ${\mathscr U}(\fg_\BC)$ be the universal enveloping
algebra of $\fg_\BC$ and ${\mathscr U}^p(\fg_\BC)$ its subspace spanned by the elements of the form
$V_1\cdots V_s$ with $V_i\in\fg_\BC$ and $s\leq p.$  We recall that there is a $\BC$-linear bijection
$\psi$ of the symmetric algebra $S(\fg_\BC)$ of $\fg_\BC$ onto ${\mathscr U}(\fg_\BC)$ which is
characterized by the property that $\psi(X^r)=X^r$ for all $X\in\fg_\BC.$ For each $\alpha\in
\textrm{Pol}_r({\mathfrak p}_\BC^*)$ we define an element $\omega(\alpha)$ of
${\mathscr U}(\fg_\BC)$ by
\begin{equation}
\omega(\alpha):=\sum \alpha_*(Y_{i_1},\cdots,Y_{i_r})\,X_{i_1}\cdots X_{i_r},
\end{equation}
where $(i_1,\cdots,i_r)$ runs over $\{1,\cdots,N\}^r.$ If $Y\in {\mathfrak p}_\BC$, then $Y^r$ as
an element of $\textrm{Pol}_r({\mathfrak p}_\BC^*)$ is defined by
$$ Y^r(u)=Y(u)^r \qquad   \textrm{for all} \ u\in {\mathfrak p}_\BC^*.$$
Hence $(Y^r)_*(u_1,\cdots,u_r)=\,Y(u_1)\cdots Y(u_r).$ According to (2.25), we see that if
$\alpha (\sum t_iY_i)=\,P(t_1,\cdots,t_N)$ for $t_i\in\BC$ with a polynomial $P$, then
\begin{equation}
\omega(\alpha)=\,\psi(P(X_1,\cdots,X_N)).
\end{equation}
Thus $\omega$ is a $\BC$-linear injection of $\textrm{Pol}({\mathfrak p}_\BC^*)$ into
${\mathscr U}(\fg_\BC)$ independent of the choice of a basis. We observe that
$\omega\big( \textrm{Pol}_r({\mathfrak p}_\BC^*)\big)=\,\psi(S_r({\mathfrak p}_\BC)).$ It is a well-known
fact that if $\alpha_1,\cdots,\alpha_m\in \textrm{Pol}_r({\mathfrak p}_\BC^*)$, then
\begin{equation}
\omega (\alpha_1\cdots\alpha_m)-\omega(\alpha_m)\cdots\omega(\alpha_1)\in {\mathscr U}^{r-1}(\fg_\BC).
\end{equation}

\vskip 0.2cm
We have a canonical pairing
\begin{equation}
\langle\, \ ,\ \,\rangle: \textrm{Pol}_r({\mathfrak p}_\BC^+) \times \textrm{Pol}_r({\mathfrak p}_\BC^-)
\lrt \BC
\end{equation}
defined by
\begin{equation}
\langle f ,g\rangle=\sum f_*({\widetilde X}_{i_1},\cdots,{\widetilde X}_{i_r}) g_*({\widetilde Y}_{i_1},\cdots,{\widetilde Y}_{i_r}),
\end{equation}
where $f_*$ (resp. $g_*$) are the unique symmetric elements of $ \textrm{Ml}_r({\mathfrak p}_\BC^+)$ (resp.\
$ \textrm{Ml}_r({\mathfrak p}_\BC^-))$, and
$\{ {\widetilde X}_1,\cdots,{\widetilde X}_{\widetilde N}\}$ and $\{ {\widetilde Y}_1,\cdots, {\widetilde Y}_{\widetilde N}\}$
are dual bases of ${\mathfrak p}_\BC^+$
and ${\mathfrak p}_\BC^-$ with respect to the Killing form $B(X,Y)=\,2(n+1)\,\textrm{tr}(XY)$, ${\widetilde N}= {{n(n+1)}\over 2},$
and $(i_1,\cdots,i_r)$ runs over $\big\{1,\cdots,{\widetilde N}\big\}^r.$

\vskip 0.3cm The adjoint representation of $K_\BC$ on ${\mathfrak p}_\BC^{\pm}$ induces the representation
of $K_\BC$ on $ \textrm{Pol}_r({\mathfrak p}_\BC^{\pm})$. Given a $K_\BC$-irreducible subspace
$Z$ of $\textrm{Pol}_r({\mathfrak p}_\BC^+),$ we can find a unique
$K_\BC$-irreducible subspace
$W$ of $\textrm{Pol}_r({\mathfrak p}_\BC^-)$ such that $\textrm{Pol}_r({\mathfrak p}_\BC^-)$ is the
direct sum of $W$ and the annihilator of $Z$. Then $Z$ and $W$ are dual with respect to the pairing (2.28).
Take bases $\{\zeta_1,\cdots,\zeta_\kappa\}$ of $Z$ and $\{ \xi_1,\cdots,\xi_\kappa\}$ of $W$ that
are dual to each other. We set
\begin{equation}
f_Z(x,y) =\,\sum_{\nu=1}^\kappa \zeta_\nu(x)\,\xi_\nu(y)\qquad (x\in {\mathfrak p}_\BC^+,
\ y\in {\mathfrak p}_\BC^-).
\end{equation}
It is easily seen that $f_Z$ belongs to $\textrm{Pol}_{2r}({\mathfrak p}_\BC)^K$ and is independent of the choice of
dual bases $\{ \zeta_\nu\}$ and $\{ \xi_\nu\}.$ Shimura \cite{Sh} proved that there exists a
canonically defined set
$\{ Z_1,\cdots,Z_n\}$ with a $K_\BC$-irreducible subspace $Z_r$ of $\textrm{Pol}_r({\mathfrak p}_\BC^+)\
(1\leq r \leq n)$ such that $f_{Z_1},\cdots,f_{Z_n}$ are algebraically independent generators of
$\textrm{Pol}({\mathfrak p}_\BC)^K$. We can identify ${\mathfrak p}_\BC^+$ with $T_n$.
We recall that $T_n$ denotes the vector space of $n\times n$ symmetric complex matrices.
We can take $Z_r$ as the subspace
of $\textrm{Pol}_r(T_n)$ spanned by the functions $f_{a;r}(Z)=\det_r(\,{}^taZa)$ for all $a\in GL(n,\BC),$
where $\det_r(x)$ denotes the determinant of the upper left $r\times r$ submatrix of $x$. For every $f\in
\textrm{Pol}({\mathfrak p}_\BC)^K$, we let $\Omega(f)$ denote the element of $\BD (\BH_n)$ represented by
$\omega (f)$. Then $\BD (\BH_n)$ is the polynomial ring $\BC[\omega (f_{Z_1}),\cdots, \omega (f_{Z_n})]$
generated by $n$ algebraically independent elements $\omega (f_{Z_1}),\cdots, \omega (f_{Z_n}).$

\end{section}

\vskip 1cm

\begin{section}{{\large\bf Invariant Differential Operators on Siegel-Jacobi Space}}
\setcounter{equation}{0}

\newcommand\Hnm{{\mathbb H}_{n,m}}
\newcommand\Cmn{\BC^{(m,n)}}
\newcommand\PO{ {{\partial}\over {\partial \Omega}} }
\newcommand\PE{ {{\partial}\over {\partial \eta}} }
\newcommand\POB{ {{\partial}\over {\partial{\overline \Omega}}} }
\newcommand\PEB{ {{\partial}\over {\partial{\overline \eta}}} }

\newcommand\pdx{ {{\partial}\over{\partial x}} }
\newcommand\pdy{ {{\partial}\over{\partial y}} }
\newcommand\pdu{ {{\partial}\over{\partial u}} }
\newcommand\pdv{ {{\partial}\over{\partial v}} }
\newcommand\PZ{ {{\partial}\over {\partial Z}} }
\newcommand\PW{ {{\partial}\over {\partial W}} }
\newcommand\PZB{ {{\partial}\over {\partial{\overline Z}}} }
\newcommand\PWB{ {{\partial}\over {\partial{\overline W}}} }
\newcommand\PX{ {{\partial\ }\over{\partial X}} }
\newcommand\PY{ {{\partial\ }\over {\partial Y}} }
\newcommand\PU{ {{\partial\ }\over{\partial U}} }
\newcommand\PV{ {{\partial\ }\over{\partial V}} }
\renewcommand\th{\theta}
\renewcommand\l{\lambda}
\renewcommand\k{\kappa}

\vskip 0.3cm The stabilizer $K^J$ of $G^J$ at $(iI_n,0)$ is given
by
\begin{equation*}
K^J=\Big\{ \big(k,(0,0;\ka)\big)\,\big|\ k\in K,\ \ka=\,^t\ka\in
\BR^{(m,m)}\,\Big\}.
\end{equation*}

\noindent Therefore $\Hnm\cong G^J/K^J$ is a homogeneous space of
$ \textit{non-reductive type}$. The Lie algebra $\fg^J$ of $G^J$
has a decomposition
\begin{equation*}
\fg^J=\fk^J+\fp^J,
\end{equation*}

\noindent where
\begin{equation*}
\fg^J=\Big\{ \big(Z,(P,Q,R)\big)\,\big|\ Z\in \fg,\ P,Q\in\BR^{(m,n)},\
R=\,^t\!R\in \BR^{(m,m)}\,\Big\},
\end{equation*}

\begin{equation*}
\fk^J=\Big\{ \big(X,(0,0,R)\big)\,\big|\ X\in \fk,\
R=\,^t\!R\in \BR^{(m,m)}\,\Big\},
\end{equation*}

\noindent
\begin{equation*}
\fp^J=\Big\{ \big(Y,(P,Q,0)\big)\,\big|\ Y\in \fp,\ P,Q\in
\BR^{(m,n)}\,\Big\}.
\end{equation*}

\noindent Thus the tangent space of the homogeneous space $\Hnm$
at $(iI_n,0)$ is identified with $\fp^J$.

\vskip 0.2cm
If $\alpha=\left( \begin{pmatrix} X_1 & \ Y_1 \\ Z_1 & -X_1 \end{pmatrix},
(P_1,Q_1,R_1)\right)$ and
$\beta=\left( \begin{pmatrix} X_2 & \ Y_2 \\ Z_2 & -X_2 \end{pmatrix},
(P_2,Q_2,R_2)\right)$ are elements of $\fg^J$, then the Lie bracket
$[\alpha,\beta]$ of $\alpha$ and $\beta$ is given by
\begin{equation}
[\alpha,\beta]=\left( \begin{pmatrix} X^* & \ Y^* \\ Z^* & -X^* \end{pmatrix},
(P^*,Q^*,R^*)\right),
\end{equation}
where
\begin{eqnarray*}
X^*&=& X_1X_2 - X_2X_1 + Y_1Z_2 - Y_2Z_1,\\
Y^*&=& X_1Y_2 - X_2Y_1 + Y_2\,^t\!X_1 - Y_1\,^t\!X_2,\\
Z^*&=& Z_1X_2 - Z_2X_1 + \,^t\!X_2 Z_1 - \,^t\!X_1 Z_2,\\
P^*&=& P_1X_2 - P_2X_1 + Q_1Z_2 - Q_2Z_1,\\
Q^*&=& P_1Y_2 - P_2Y_1 + Q_2\,^t\!X_1 - Q_1\,^t\!X_2,\\
R^*&=& P_1\,^t\!Q_2 - P_2\,^t\!Q_1 + Q_2\,^t\!P_1 - Q_1\,^t\!P_2
\end{eqnarray*}

\begin{lemma} \begin{equation*}
[\fk^J,\fk^J]\subset \fk^J,\quad [\fk^J,\fp^J]\subset \fp^J.
\end{equation*}
\end{lemma}
\vskip 0.2cm \noindent {\it Proof.} The proof follows immediately from Formula (3.1).
$\hfill\square$

\vskip 0.3cm

\begin{lemma}
Let
$$k^J = \left( \begin{pmatrix} A & -B \\ B & \ A \end{pmatrix},
(0,0,\kappa)\right)\in K^J $$
with $ \begin{pmatrix} A & -B \\ B & \ A \end{pmatrix}\in K,
\ \ka=\,^t\ka \in\BR^{(m,m)}$
and
$$\alpha=\left( \begin{pmatrix} X & \ Y\\ Y & -X \end{pmatrix},(P,Q,0)\right)
\in \fp^J$$
with $ X=\,^tX,\ Y=\,^tY\in \BR^{(n,n)},\ P,Q\in\BR^{(m,n)}.$
Then the adjoint action of $K^J$ on $\fp^J$ is given by
\begin{equation}
\textrm{Ad}(k^J)\alpha = \left( \begin{pmatrix} X_* & \ Y_*\\ Y_* & -X_* \end{pmatrix},(P_*,Q_*,0)\right),
\end{equation}
where
\begin{eqnarray}
X_* &=& AX\,^t\!A - \big( BX\,^t\!B + BY\,^t\!A + AY\,^t\!B\big),\\
Y_* &=& \big( AX\,^t\!B + AY\,^t\!A + BX\,^t\!A\big) -BY\,^t\!B,\\
P_* &=& P\,\,^t\!A - Q\,\,^t\!B,\\
Q_* &=& P\,\,^t\!B + Q\,\,^t\!A.
\end{eqnarray}
\end{lemma}
\vskip 0.2cm
\noindent {\it Proof.} We leave the proof to the reader. $\hfill\square$

\vskip 0.3cm
We recall that $T_n$ denotes the vector space of all $n\times n$ symmetric complex matrices.
For brevity, we put $T_{n,m}:=T_n\times \BC^{(m,n)}.$
We define the real linear map $\Phi:\fp^J\lrt T_{n,m}$ by
\begin{equation}
\Phi \left( \begin{pmatrix} X & \ Y\\ Y & -X \end{pmatrix},(P,Q,0)\right)=
\big( X\,+\,i\,Y,\,P\,+\,i\,Q \big),
\end{equation}
where $\begin{pmatrix} X & \ Y\\ Y & -X\end{pmatrix}\in \fp$ and $P,Q\in \BR^{(m,n)}.$

\vskip 0.3cm Let $S(n,\BR)$ denote the additive group consisting of all $n\times n$
real symmetric matrices. Now we define the isomorphism $\theta : K^J \lrt U(n)\times S(n,\BR)$
by
\begin{equation}
\theta (h,(0,0,\kappa))= (\delta(h), \kappa),\quad h\in K, \ \kappa\in S(n,\BR),
\end{equation}
where $\delta :K\lrt U(n)$ is the map defined by (2.5).
Identifying $\BR^{(m,n)}\times \BR^{(m,n)}$ with
$\BC^{(m,n)}$, we can identify $\fp^J$ with $T_n\times
\BC^{(m,n)}$.

\vskip 0.3cm
\begin{theorem}
The adjoint representation of $K^J$ on $\fp^J$ is compatible with the {\it natural\ action} of
$U(n)\times S(n,\BR)$ on $T_{n,m}$ defined by
\begin{equation}
(h,\kappa)\cdot (\omega,z):= (h\,\omega \,^th,\,z\,^th),\qquad
h\in U(n),\ \kappa\in S(n,\BR), \ (\omega,z)\in T_{n,m}
\end{equation}
through the maps $
\Phi$ and $\theta$. Precisely, if $k^J\in K^J$ and $\alpha\in \fp^J$, then we have the
following equality
\begin{equation}
\Phi \big( Ad\big( k^J\,\big)\alpha\Big) = \theta \big( k^J\, \big)\cdot \Phi (\alpha).
\end{equation}
Here we
regard the complex vector space $T_{n,m}$ as a real vector space.
\end{theorem}

\vskip 0.3cm \noindent
{\it Proof.} Let
$$k^J = \left( \begin{pmatrix} A & -B \\ B & \ A \end{pmatrix},
(0,0,\kappa)\right)\in K^J $$
with $ \begin{pmatrix} A & -B \\ B & \ A \end{pmatrix}\in K,
\ \ka=\,^t\ka \in\BR^{(m,m)}$
and
$$\alpha=\left( \begin{pmatrix} X & \ Y\\ Y & -X \end{pmatrix},(P,Q,0)\right)
\in \fp^J$$
with $ X=\,^tX,\ Y=\,^tY\in \BR^{(n,n)},\ P,Q\in\BR^{(m,n)}.$
Then we have

\begin{eqnarray*}
\theta \big( k^J\, \big)\cdot \Phi (\alpha)
&=& \big( A\,+\,i\,B,\,\kappa\big)\cdot \big(  X\,+\,i\,Y,\,P\,+\,i\,Q\big)\\
&=& \big( (A+iB)(X+iY)\,^t(A+iB),\,(P+iQ)\,^t\!(A+iB)\big)\\
&=& \big( X_*\,+\,i\,Y_*,\,P_*\,+\,i\,Q_*\big)\\
&=& \Phi  \left( \begin{pmatrix} X_* & \ Y_*\\ Y_* & -X_* \end{pmatrix},(P_*,Q_*,0)\right)\\
&=& \Phi \big( Ad\big( k^J\,\big)\alpha\Big)\qquad (by \ Lemma\ 3.2),
\end{eqnarray*}
where $X_*,Y_*,Z_*$ and $Q_*$ are given by the formulas (3.3),\,(3.4),\,(3.5) and (3.6)
respectively.
$\hfill\square$

\newcommand\bw{d{\overline W}}
\newcommand\bz{d{\overline Z}}
\newcommand\bo{d{\overline \Omega}}

\vskip 0.3cm We now study the algebra $\BD(\BH_{n,m})$ of all
differential operators on $\BH_{n,m}$ invariant under the {\it natural action}
(1.2) of $G^J$. The action (3.9) induces the action of $U(n)$ on
the polynomial algebra $\text{Pol}_{n,m}:=\,\text{Pol}\,(T_{n,m}).$
We denote by $\text{Pol}_{n,m}^{U(n)}$ the subalgebra of
$\text{Pol}_{n,m}$ consisting of all $U(n)$-invariants. Similarly
the action (3.2) of $K$ induces the action of $K$ on the
polynomial algebra $ \textrm{Pol}\big(\fp^J\big)$. We see that
through the identification of $\fp^J$ with $T_{n,m}$, the algebra
$ \textrm{Pol}\big(\fp^J\big)$ is isomorphic to
$\text{Pol}_{n,m}.$ The following $U(n)$-invariant inner product
$(\,\,,\,)_*$ of the complex vector space $T_{n,m}$ defined by
\begin{equation*}
\big((\om,z),(\om',z')\big)_*= \textrm{tr}\big(\om{\overline
{\om'}}\,\big)+ \textrm{tr}\big(z\,^t{\overline {z'}}\,\big),\quad
(\om,z),\,(\om',z')\in T_{n,m}
\end{equation*}

\noindent gives a canonical isomorphism
\begin{equation*}
T_{n,m}\cong\,T_{n,m}^*,\quad (\om,z)\mapsto f_{\om,z},\quad
(\om,z)\in T_{n,m},
\end{equation*}

\noindent where $f_{\om,z}$ is the linear functional on $T_{n,m}$
defined by
\begin{equation*}
f_{\om,z}\big((\om',z'\,)\big)=\big((\om',z'),(\om,z)\big)_*,\quad
(\om',z'\,)\in T_{n,m}.
\end{equation*}

\noindent According to Helgason (\cite{He2}, p.\,287), one gets a canonical linear bijection
of $S(T_{n,m})^{U(n)}$ onto $\BD(\BH_{n,m})$. Identifying
$T_{n,m}$ with $T_{n,m}^*$ by the above isomorphism, one gets a
natural linear bijection
$$\Theta_{n,m}:\,\text{Pol}^{U(n)}_{n,m}\lrt \BD(\BH_{n,m})$$
of $\text{Pol}^{U(n)}_{n,m}$ onto $\BD(\BH_{n,m}).$ The map
$\Theta_{n,m}$ is described explicitly as follows. We put
$N_{\star}=n(n+1)+2mn$. Let $\big\{ \eta_{\alpha}\,|\ 1\leq \alpha
\leq N_{\star}\, \big\}$ be a basis of $\fp^J$. If $P\in
\textrm{Pol}\big(\fp^J\big)^K=\mathrm{Pol}_{n,m}^{U(n)}$, then
\begin{equation}
\Big(\Theta_{n,m} (P)f\Big)(gK^J)=\left[ P\left( {{\partial}\over
{\partial t_{\al}}}\right)f\left(g\,\text{exp}\,
\left(\sum_{\al=1}^{N_{\star}} t_{\al}\eta_{\al}\right)
K^J\right)\right]_{(t_{\al})=0},
\end{equation}

\noindent where $f\in C^{\infty}({\mathbb H}_{n,m})$. In general,
it is hard to express $\Theta_{n,m}(P)$ explicitly for a polynomial
$P\in \textrm{Pol}\big(\fp^J\big)^K$. We refer to \cite{He2}, p. 287.

\vskip 0.2cm  We present the following $\textit{basic}$
$U(n)$-invariant polynomials
in $\text{Pol}_{n,m}^{U(n)}$.
\begin{eqnarray}
&& q_j(\om,z)=\,\text{tr}\big((\om\,{\overline
\om})^{j+1}\big),\quad 0\leq j\leq n-1,\\
&&  \alpha_{kp}^{(j)}(\om,z)=
\,\text{Re}\,\big( z\,({\overline \om}\om)^j \,^t{\overline z} \big)_{kp},
\quad 0\leq j\leq n-1,\ 1\leq k\leq p\leq m,\\
&& \beta_{lq}^{(j)}(\om,z)=
\,\text{Im}\,\big( z\,({\overline \om}\om)^j \,^t{\overline z} \big)_{lq},
\quad 0\leq j\leq n-1
,\ 1\leq l< q\leq m,\\
&& f_{kp}^{(j)}(\om,z)= \,\text{Re}\,(z\,({\overline \om}\om)^j \,
{\overline
\om}\,^t\!z)_{kp}, \quad   0\leq j\leq n-1,\ 1\leq k\leq p\leq m,\\
&& g_{kp}^{(j)}(\om,z)= \,\text{Im}\,(z\,
({\overline \om}\om)^j \,{\overline
\om}\,^t\!z\,)_{kp},\quad  0\leq j\leq n-1,\ 1\leq k\leq p\leq m,
\end{eqnarray}
\noindent where $\om\in T_n$ and $z\in \BC^{(m,n)}$.

\vskip 0.3cm We present some interesting $U(n)$-invariants.
For an $m\times m$ matrix $S$, we define the
following invariant polynomials in $\text{Pol}_{n,m}^{U(n)}$:
\begin{eqnarray}
&&
m_{j;S}^{(1)}(\om,z)=\,\textrm{Re}\,\Big(\text{tr}\big(\om{\overline
\om}+ \,^tzS{\overline
z}\,\big)^j\,\Big),\quad 1\leq j\leq n,\\
&&
m_{j;S}^{(2)}(\om,z)=\,\textrm{Im}\,\Big(\text{tr}\big(\om{\overline
\om}+ \,^tzS{\overline
z}\,\big)^j\,\Big),\quad 1\leq j\leq n,\\
&&  q_{k;S}^{(1)}(\om,z)=\,\textrm{Re}\,\Big( \textrm{tr}\big( (
\,^tz\,S\,{\overline
z})^k\big) \Big),\quad 1\leq k\leq m,  \\
&&  q_{k;S}^{(2)}(\om,z)=\,\textrm{Im}\,\Big( \textrm{tr}\big( (
\,^tz\,S\,{\overline
z})^k\big) \Big),\quad 1\leq k\leq m,  \\
&& \theta_{i,k,j;S}^{(1)}(\om,z) =\,\textrm{Re}\,\Big(
\textrm{tr}\big( (\om {\overline \om})^i\,(\,^tz\,S\,{\overline
z})^k\,(\om {\overline \om}+\,^tz\,S\,{\overline
z}\,)^j\,\big)\Big),\\
&& \theta_{i,k,j;S}^{(2)}(\om,z) =\,\textrm{Im}\,\Big(
\textrm{tr}\big( (\om {\overline \om})^i\,(\,^tz\,S\,{\overline
z})^k\,(\om {\overline \om}+\,^tz\,S\,{\overline
z}\,)^j\,\big)\Big),
\end{eqnarray}
\noindent where $1\leq i,j\leq n$ and $1\leq k\leq m$.

\vskip 0.2cm We define the following $U(n)$-invariant polynomials in
$\text{Pol}_{n,m}^{U(n)}$.
\begin{eqnarray}
&&  r_{jk}^{(1)}(\om,z)= \,\textrm{Re}\,\Big( \textrm{det}\big(
(\om {\overline \om})^j\,(\,^tz{\overline z})^k\,\big)\Big),
\quad 1\leq j\leq n,\ 1\leq k\leq m,\\
&&  r_{jk}^{(2)}(\om,z)= \,\textrm{Im}\,\Big( \det\big(
(\om {\overline \om})^j\,(\,^tz{\overline z})^k\,\big)\Big), \quad
1\leq j\leq n,\ 1\leq k\leq m.
\end{eqnarray}

 \vskip 0.355cm We propose the following natural problems.

\vskip 0.2cm \noindent $ \textbf{Problem 1.}$ Find a complete list of explicit generators
of $\text{Pol}_{n,m}^{U(n)}$.

\vskip 0.32cm \noindent $ \textbf{Problem 2.}$ Find all the relations among a set of generators of $\text{Pol}_{n,m}^{U(n)}$.

\vskip 0.32cm \noindent $ \textbf{Problem 3.}$ Find an easy or effective way to
express the images of the above invariant polynomials or generators of $\text{Pol}_{n,m}^{U(n)}$ under the
Helgason map $\Theta_{n,m}$ explicitly.

\vskip 0.32cm \noindent $ \textbf{Problem 4.}$ Decompose $\text{Pol}_{n,m}$ into
$U(n)$-irreducibles.

\vskip 0.32cm \noindent $ \textbf{Problem 5.}$ Find a complete list of
 explicit generators of the algebra $\BD(\BH_{n,m})$. Or
construct explicit $G^J$-invariant differential operators on $\BH_{n,m}.$

\vskip 0.32cm \noindent $ \textbf{Problem 6.}$ Find all the relations among a set of generators of $\BD(\BH_{n,m})$.

\vskip 0.32cm \noindent $ \textbf{Problem 7.}$ Is $\text{Pol}_{n,m}^{U(n)}$ finitely generated ? Is $\BD(\BH_{n,m})$
finitely generated ?


\vskip 0.5cm
Quite recently Minoru Itoh \cite{It} solved Problem 1 and Problem 7.

\begin{theorem} $\textrm{Pol}_{n,m}^{U(n)}$ is generated by
$$q_j(\om,z),\ \alpha_{kp}^{(j)}(\om,z),\ \beta_{lq}^{(j)}(\om,z)
,\ f_{kp}^{(j)}(\om,z)\ \text{and}\ g_{kp}^{(j)}(\om,z),$$
where $0\leq j\leq n-1,\ 1\leq k\leq p\leq m \ \textrm{and}\ \ 1\leq l< q\leq m$.
\end{theorem}

\end{section}

\newcommand\Hnm{{\mathbb H}_{n,m}}
\newcommand\Cmn{\BC^{(m,n)}}
\newcommand\PO{ {{\partial}\over {\partial \Omega}} }
\newcommand\PE{ {{\partial}\over {\partial \eta}} }
\newcommand\POB{ {{\partial}\over {\partial{\overline \Omega}}} }
\newcommand\PEB{ {{\partial}\over {\partial{\overline \eta}}} }

\newcommand\pdx{ {{\partial}\over{\partial x}} }
\newcommand\pdy{ {{\partial}\over{\partial y}} }
\newcommand\pdu{ {{\partial}\over{\partial u}} }
\newcommand\pdv{ {{\partial}\over{\partial v}} }
\newcommand\PZ{ {{\partial}\over {\partial Z}} }
\newcommand\PW{ {{\partial}\over {\partial W}} }
\newcommand\PZB{ {{\partial}\over {\partial{\overline Z}}} }
\newcommand\PWB{ {{\partial}\over {\partial{\overline W}}} }
\newcommand\PX{ {{\partial\ }\over{\partial X}} }
\newcommand\PY{ {{\partial\ }\over {\partial Y}} }
\newcommand\PU{ {{\partial\ }\over{\partial U}} }
\newcommand\PV{ {{\partial\ }\over{\partial V}} }

\newcommand\bw{d{\overline W}}
\newcommand\bz{d{\overline Z}}
\newcommand\bo{d{\overline \Omega}}

\renewcommand\th{\theta}
\renewcommand\l{\lambda}
\renewcommand\k{\kappa}

\begin{section}{{\large\bf Examples of Explicit $G^J$-Invariant Differential Operators }}
\setcounter{equation}{0}
\vskip 0.3cm In this section we give examples of explicit $G^J$-invariant differential operators on
the Siegel-Jacobi space and the Siegel-Jacobi disk.
\vskip 0.3cm For $g=\big(M,(\la,\mu;\kappa)\big)\in G^J$ with $M=
\begin{pmatrix} A & B \\ C & D \end{pmatrix}\in Sp(n,\BR)$ and $(\Omega,Z)\in \BH_{n,m},$
we set
\begin{eqnarray*}
\Omega_*&=&\,M\!\cdot\!\Omega=\,X_*+i\,Y_*,\quad X_*,Y_*\ \textrm{real},\\
Z_*&=&\,(Z+\la \Omega+\mu)(C\Omega+D)^{-1}=\,U_*+ i\,V_*,\quad U_*,V_*\ \textrm{real}.
\end{eqnarray*}
For a coordinate
$(\Om,Z)\in \BH_{n,m}$ with $\Om=(\omega_{\mu\nu})$ and
$Z=(z_{kl})$, we put $d\Om,\,d{\overline \Om},\,\PO,\,\POB$ as
before and set
\begin{eqnarray*}
Z\,&=&U\,+\,iV,\quad\ \ U\,=\,(u_{kl}),\quad\ \ V\,=\,(v_{kl})\ \
\text{real},\\
dZ\,&=&\,(dz_{kl}),\quad\ \ d{\overline Z}=(d{\overline z}_{kl}),
\end{eqnarray*}
$$\PZ=\begin{pmatrix} {\partial}\over{\partial z_{11}} & \hdots &
 {\partial}\over{\partial z_{m1}} \\
\vdots&\ddots&\vdots\\
 {\partial}\over{\partial z_{1n}} &\hdots & {\partial}\over
{\partial z_{mn}} \end{pmatrix},\quad \PZB=\begin{pmatrix}
{\partial}\over{\partial {\overline z}_{11} }   &
\hdots&{ {\partial}\over{\partial {\overline z}_{m1} }  }\\
\vdots&\ddots&\vdots\\
{ {\partial}\over{\partial{\overline z}_{1n} }  }&\hdots &
 {\partial}\over{\partial{\overline z}_{mn} }  \end{pmatrix}.$$
Then we can show that
\begin{eqnarray}
d\Omega_*&=&\,{}^t(C\Om+D)^{-1}\,d\Om (C\Omega+D)^{-1},\\
dZ_* &=&\,dZ(C\Om+D)^{-1}\,\\
 & & \ \ +\,\left\{ \la-(Z+\la \Omega+\mu)(C\Om+D)^{-1}C\right\}d\Om (C\Om+D)^{-1},\nonumber\\
{{\partial\ \,}\over {\partial{\Om_*}} }&=&\, (C\Om+D)\,{}^{{}^{{}^{{}^\text{\scriptsize $t$}}}}\!\!\!
\left\{ (C\Om+D) {{\partial\ }\over {\partial \Om}} \right\}\\
& &\ \ + (C\Om+D)\,{}^{{}^{{}^{{}^\text{\scriptsize $t$}}}}\!\!\!\left\{
\big( C\,{}^t\!Z+C\,{}^t\mu-D\,{}^t\!\lambda\big) \,{}^{{}^{{}^{{}^\text{\scriptsize $t$}}}}\!\!\!\left(
\PZ\right)\right\}\nonumber
\end{eqnarray}
and
\begin{equation}
{{\partial\ \,}\over {\partial{Z_*}} }=\,(C\Om+D)\PZ.
\end{equation}
\vskip 0.2cm From \cite[p.\,33]{M2} or \cite[p.\,128]{Si2}, we know that
\begin{equation}
Y_*=\,{}^t(C{\overline\Om}+D)^{-1} Y (C\Om+D)^{-1}=\,{}^t(C\Om+D)^{-1}Y (C{\overline\Om}+D)^{-1}.
\end{equation}

Using Formulas (4.1),\,(4.2) and (4.5), the author \cite{YJH7} proved that for any two positive real numbers
$A$ and $B$,
\begin{eqnarray*}
ds_{n,m;A,B}^2&=&\,A\, \tra\Big( Y^{-1}d\Om\,Y^{-1}d{\overline
\Om}\Big) \nonumber \\ && \ \ + \,B\,\bigg\{ \tra\Big(
Y^{-1}\,^tV\,V\,Y^{-1}d\Om\,Y^{-1} d{\overline \Om} \Big)
 +\,\tra\Big( Y^{-1}\,^t(dZ)\,\bz\Big)\\
&&\quad\quad -\tra\Big( V\,Y^{-1}d\Om\,Y^{-1}\,^t(\bz)\Big)\,
-\,\tra\Big( V\,Y^{-1}d{\overline \Om}\, Y^{-1}\,^t(dZ)\,\Big)
\bigg\} \nonumber
\end{eqnarray*}
is a Riemannian metric on $\BH_{n,m}$ which is invariant under the action (1.2) of $G^J.$

\vskip 0.3cm The following lemma is very useful for computing the invariant differential operators.
H. Maass \cite{M} observed the following useful fact. \vskip 0.3cm
\begin{lemma} (a) Let $A$ be an $m\times n$ matrix and
$B$ an $n\times l$ matrix. Assume that the entries of $A$ commute
with the entries of $B$. Then ${}^t(AB)=\,{ }^tB\,\,{ }^tA.$
\vskip 0.1cm \noindent (b) Let $A,\,B$ and $C$ be a $k\times l$,
an $n\times m$ and an $m\times l$ matrix respectively. Assume that
the entries of $A$ commute with the entries of $B$. Then
\begin{eqnarray*}
{ }^t(A\,\,{ }^t(BC))=\,B\,\,{ }^t(A\,\,^tC).
\end{eqnarray*}
\end{lemma}
\vskip 0.2cm \noindent {\it Proof.} The proof follows immediately
from the direct computation. \hfill$\Box$

\vskip 0.3cm
Using Formulas (4.3),\,(4.4),\,(4.5) and Lemma 4.1, the author \cite{YJH7} proved that the following differential operators
${\mathbb M}_1$ and ${\mathbb M}_2$
on $\BH_{n,m}$ defined by
\begin{equation}
{\mathbb M}_1=\,\textrm{tr}\left(\, Y\,\PZ\,
{}^{{}^{{}^{{}^\text{\scriptsize $t$}}}}\!\!\!
\left(
\PZB\right)\,\right)
\end{equation}
and
\begin{eqnarray}
{\mathbb M}_2\,&=&  \tra\left(\,Y\,
{}^{{}^{{}^{{}^\text{\scriptsize $t$}}}}\!\!\!
\left(Y\POB\right)\PO\,\right)\, +\,\tra\left(\,VY^{-1}\,^tV\,
{}^{{}^{{}^{{}^\text{\scriptsize $t$}}}}\!\!\!
\left(Y\PZB\right)\,\PZ\,\right)\\
& &\ \
+\,\tra\left(V\,
{}^{{}^{{}^{{}^\text{\scriptsize $t$}}}}\!\!\!
\left(Y\POB\right)\PZ\,\right)
+\,\tra\left(\,^tV\,
{}^{{}^{{}^{{}^\text{\scriptsize $t$}}}}\!\!\!
\left(Y\PZB\right)\PO\,\right)\nonumber
\end{eqnarray}
are invariant under the action (1.2) of $G^J.$  The author \cite{YJH7} proved that
for any two positive real numbers
$A$ and $B$, the following differential operator

\begin{equation}
\Delta_{n,m;A,B}=\,{\frac 4A}\,{\mathbb M}_2 + {\frac 4B}
{\mathbb M}_1
\end{equation}
is the Laplacian of the $G^J$-invariant Riemannian metric $ds_{n,m;A,B}^2.$

\vskip 0.5cm
\begin{proposition} The following differential
operator ${\mathbb K}$ on $\Hnm$ of degree $2n$ defined by
\begin{equation}
{\mathbb K}=\,\det(Y)\,\det\left( \PZ {}^{{}^{{}^{{}^\text{\scriptsize $t$}}}}\!\!\!\left(
\PZB\right)\right)
\end{equation}
is invariant under the action (1.2) of $G^J$.
\end{proposition}

\vskip 0.2cm\noindent {\it Proof.} Let ${\mathbb K}_{M,(\la,\mu;\kappa)}$ denote the image of ${\mathbb K}$ under the transformation
\begin{equation*}
(\Omega,Z)\longmapsto \big((M\!\cdot\! \Omega,(Z+\la \Om+\mu)(C\Om+D)^{-1}\big)
\end{equation*}
with $M=\begin{pmatrix} A & B\\ C& D \end{pmatrix}\in Sp(n,\BR)$ and
$(\la,\mu;\kappa)\in H_\BR^{(n,m)}.$
If $f$ is a $C^\infty$ function on $\Hnm,$ using (4.4),\,(4.5) and Lemma 4.1, we have
\begin{eqnarray*}
{\mathbb K}_{M,(\la,\mu;\kappa)}f&=&\,\det(Y)\,|\det (C\Omega+D)|^{-2}\,\det\left[ (C\Omega+D)\PZ
{}^{{}^{{}^{{}^\text{\scriptsize $t$}}}}\!\!\!\left\{ (C{\overline\Om}+D)
{{\partial f}\over{\partial {\overline Z}}}   \right\}\right]\\
&=&\,\det(Y)\,|\det (C\Omega+D)|^{-2}\,\det\left[ (C\Omega+D)
{}^{{}^{{}^{{}^\text{\scriptsize $t$}}}}\!\!\!\left\{ (C{\overline\Om}+D)
{}^{{}^{{}^{{}^\text{\scriptsize $t$}}}}\!\!\!\left( \PZ {}^{{}^{{}^{{}^\text{\scriptsize $t$}}}}\!\!\!
\left( {{\partial f}\over{\partial {\overline Z}}} \right)\right)  \right\}\right]\\
&=&\,\det(Y)\,|\det (C\Omega+D)|^{-2}\,\det\left[ (C\Omega+D)\PZ
{}^{{}^{{}^{{}^\text{\scriptsize $t$}}}}\!\!\!\left( {{\partial f}\over{\partial {\overline Z}}} \right)
\,{}^t(C{\overline\Om}+D)\right]\\
&=&\,\det(Y)\,\det \left( \PZ {}^{{}^{{}^{{}^\text{\scriptsize $t$}}}}\!\!\!\left(
{{\partial f}\over{\partial {\overline Z}}} \right) \right)\\
&=&\,{\mathbb K}f.
\end{eqnarray*}
Since $M\in Sp(n,\BR)$ and $(\la,\mu;\kappa)\in H_\BR^{(n,m)}$ are arbitrary, $\mathbb K$ is invariant under
the action (1.2) of $G^J$.
\hfill$\square$

\vskip 0.5cm
\begin{proposition} The following
matrix-valued differential operator ${\mathbb T}$ on
$\BH_{n,m}$ defined by
\begin{equation}
{\mathbb T}=\,
{}^{{}^{{}^{{}^\text{\scriptsize $t$}}}}\!\!\!
\left( \PZB\right) Y \PZ
\end{equation}

\noindent is invariant under the action (1.2) of $G^J$.
\end{proposition}
\vskip 0.2cm\noindent
{\it Proof.} Let ${\mathbb T}_{M,(\la,\mu;\kappa)}$ denote the image of ${\mathbb K}$ under the transformation
\begin{equation*}
(\Omega,Z)\longmapsto \big((M\!\cdot\! \Omega,(Z+\la \Om+\mu)(C\Om+D)^{-1}\big)
\end{equation*}
with $M=\begin{pmatrix} A & B\\ C& D \end{pmatrix}\in Sp(n,\BR)$ and
$(\la,\mu;\kappa)\in H_\BR^{(n,m)}.$
If $f$ is a $C^\infty$ function on $\Hnm,$ according to (4.4), (4.5) and Lemma 4.1, we have
\begin{eqnarray*}
{\mathbb T}_{M,(\la,\mu;\kappa)}f&=&\,{}^{{}^{{}^{{}^\text{\scriptsize $t$}}}}\!\!\!\left(
(C{\overline \Omega}+D){{\partial\ }\over{\partial {\overline Z}}}   \right)\,{}^t\!(C{\overline\Omega}+D)^{-1}
Y (C\Omega+D)^{-1} (C\Omega+D)
{{\partial f}\over{\partial Z}}\\
&=&\,
{}^{{}^{{}^{{}^\text{\scriptsize $t$}}}}\!\!\!
\left( \PZB\right) Y {{\partial f}\over{\partial Z}}\\
&=&\,{\mathbb T}f.\\
\end{eqnarray*}
Since $M\in Sp(n,\BR)$ and $(\la,\mu;\kappa)\in H_\BR^{(n,m)}$ are arbitrary, $\mathbb T$ is invariant under
the action (1.2) of $G^J$.
\hfill$\square$

\vskip 0.2cm\noindent
\begin{corollary} Each $(k,l)$-entry ${\mathbb T}_{kl}$ of ${\mathbb
T}$ given by
\begin{equation}
{\mathbb T}_{kl}=\sum_{i,j=1}^n
\,y_{ij}\,{{\partial^2\ \ \ \ }\over{\partial {\overline z}_{ki}\partial
z_{lj}} },\quad 1\leq k,l\leq m
\end{equation}

\noindent is an element of $\BD\big(\BH_{n,m}\big)$.
\end{corollary}
\vskip 0.2cm\noindent
{\it Proof.} It follows immediately from Proposition 4.2. \hfill$\square$


\newcommand\OW{\overline{W}}
\newcommand\OP{\overline{P}}
\newcommand\OQ{\overline{Q}}
\newcommand\OZ{\overline{Z}}
\newcommand\Dn{{\mathbb D}_n}
\newcommand\Dnm{{\mathbb D}_{n,m}}
\newcommand\Hn{{\mathbb H}_n}
\newcommand\SJ{{\mathbb H}_n\times {\mathbb C}^{(m,n)}}
\newcommand\DC{{\mathbb D}_n\times {\mathbb C}^{(m,n)}}
\newcommand\ot{\overline\eta}


\vskip 0.3cm Now we consider invariant differential operators on the Siegel-Jacobi disk.
Let
\begin{equation*}
\BD_n=\left\{ W\in \BC^{(n,n)}\,|\ W=\,{}^tW,\ I_n-\OW W >
0\,\right\}
\end{equation*}
be the generalized unit disk.

\vskip 0.2cm For brevity, we write $\Dnm:=\DC.$ For a coordinate
$(W,\eta)\in\Dnm$ with $W=(w_{\mu\nu})\in {\mathbb D}_n$ and
$\eta=(\eta_{kl})\in \Cmn,$ we put
\begin{eqnarray*}
dW\,&=&\,(dw_{\mu\nu}),\quad\ \ d{\overline W}\,=\,(d{\overline w}_{\mu\nu}),\\
d\eta\,&=&\,(d\eta_{kl}),\quad\ \
d{\overline\eta}\,=\,(d{\overline\eta}_{kl})
\end{eqnarray*}
and
\begin{eqnarray*}
\PW\,=\,\left(\, { {1+\delta_{\mu\nu}} \over 2}\, {
{\partial}\over {\partial w_{\mu\nu}} } \,\right),\quad
\PWB\,=\,\left(\, { {1+\delta_{\mu\nu}}\over 2} \, {
{\partial}\over {\partial {\overline w}_{\mu\nu} }  } \,\right),
\end{eqnarray*}

$$\PE=\begin{pmatrix} {\partial}\over{\partial \eta_{11}} & \hdots &
 {\partial}\over{\partial \eta_{m1}} \\
\vdots&\ddots&\vdots\\
 {\partial}\over{\partial \eta_{1n}} &\hdots & {\partial}\over
{\partial \eta_{mn}} \end{pmatrix},\quad \PEB=\begin{pmatrix}
{\partial}\over{\partial {\overline \eta}_{11} }   &
\hdots&{ {\partial}\over{\partial {\overline \eta}_{m1} }  }\\
\vdots&\ddots&\vdots\\
{ {\partial}\over{\partial{\overline \eta}_{1n} }  }&\hdots &
 {\partial}\over{\partial{\overline \eta}_{mn} }  \end{pmatrix}.$$

\vskip 0.2cm  We can identify an element $g=(M,(\la,\mu;\kappa))$
of $G^J,\ M=\begin{pmatrix} A&B\\
C&D\end{pmatrix}\in Sp(n,\BR)$ with the element
\begin{equation*}
\begin{pmatrix} A & 0 & B & A\,^t\mu-B\,^t\la  \\ \la & I_m & \mu
& \kappa \\ C & 0 & D & C\,^t\mu-D\,^t\la \\ 0 & 0 & 0 & I_m
\end{pmatrix}
\end{equation*}
of $Sp(m+n,\BR).$ \vskip 0.3cm We set
\begin{equation*}
T_*={1\over {\sqrt 2}}\,
\begin{pmatrix} I_{m+n} & I_{m+n}\\ iI_{m+n} & -iI_{m+n}
\end{pmatrix}.
\end{equation*}
We now consider the group $G_*^J$ defined by
\begin{equation*}
G_*^J:=T_*^{-1}G^JT_*.
\end{equation*}
If $g=(M,(\la,\mu;\kappa))\in G^J$ with $M=\begin{pmatrix} A&B\\
C&D\end{pmatrix}\in Sp(n,\BR)$, then $T_*^{-1}gT_*$ is given by
\begin{equation}
T_*^{-1}gT_*=
\begin{pmatrix} P_* & Q_*\\ {\overline Q}_* & {\overline P}_*
\end{pmatrix},
\end{equation}
where
\begin{equation*}
P_*=
\begin{pmatrix} P & {\frac 12} \left\{ Q\,\,{}^t(\la+i\mu)-P\,\,{}^t(\la-i\mu)\right\}\\
{\frac 12} (\la+i\mu) & I_h+i{\frac \kappa 2}
\end{pmatrix},
\end{equation*}

\begin{equation*}
Q_*=
\begin{pmatrix} Q & {\frac 12} \left\{ P\,\,{}^t(\la-i\mu)-Q\,\,{}^t(\la+i\mu)\right\}\\
{\frac 12} (\la-i\mu) & -i{\frac \kappa 2}
\end{pmatrix},
\end{equation*}
and $P,\,Q$ are given by the formulas
\begin{equation}
P= {\frac 12}\,\left\{ (A+D)+\,i\,(B-C)\right\}
\end{equation}
and
\begin{equation}
 Q={\frac
12}\,\left\{ (A-D)-\,i\,(B+C)\right\}.
\end{equation}
From now on, we write
\begin{equation*}
\left(\begin{pmatrix} P & Q\\ {\overline Q} & {\overline P}
\end{pmatrix},\left( {\frac 12}(\la+i\mu),\,{\frac 12}(\la-i\mu);\,-i{\kappa\over 2}\right)\right):=
\begin{pmatrix} P_* & Q_*\\ {\overline Q}_* & {\overline P}_*
\end{pmatrix}.
\end{equation*}
In other words, we have the relation
\begin{equation*}
T_*^{-1}\left( \begin{pmatrix} A&B\\
C&D\end{pmatrix},(\la,\mu;\kappa)
\right)T_*=  \left(\begin{pmatrix} P & Q\\
{\overline Q} & {\overline P}
\end{pmatrix},\left(
{\frac 12}(\la+i\mu),\,{\frac 12}(\la-i\mu);\,-i{\kappa\over 2}
\right)\right).
\end{equation*}
Let
\begin{equation*}
H_{\BC}^{(n,m)}:=\left\{ (\xi,\eta\,;\zeta)\,|\
\xi,\eta\in\BC^{(m,n)},\ \zeta\in\BC^{(m,m)},\
\zeta+\eta\,{}^t\xi\ \textrm{symmetric}\,\right\}
\end{equation*}
be the complex Heisenberg group endowed with the following
multiplication
\begin{equation*}
(\xi,\eta\,;\zeta)\circ
(\xi',\eta';\zeta'):=(\xi+\xi',\eta+\eta'\,;\zeta+\zeta'+
\xi\,{}^t\eta'-\eta\,{}^t\xi')).
\end{equation*}
We define the semidirect product
\begin{equation*}
SL(2n,\BC)\ltimes H_{\BC}^{(n,m)}
\end{equation*}
endowed with the following multiplication
\begin{eqnarray*}
& & \left( \begin{pmatrix} P & Q\\ R & S
\end{pmatrix}, (\xi,\eta\,;\zeta)\right)\cdot \left( \begin{pmatrix} P' & Q'\\
R' & S'
\end{pmatrix}, (\xi',\eta';\zeta')\right)\\
&=& \left( \begin{pmatrix} P & Q\\ R & S
\end{pmatrix}\,\begin{pmatrix} P' & Q'\\ R' & S'
\end{pmatrix},\,({\tilde \xi}+\xi',{\tilde
\eta}+\eta';\zeta+\zeta'+{\tilde \xi}\,{}^t\eta'-{\tilde
\eta}\,{}^t\xi')  \right),
\end{eqnarray*}
where ${\tilde\xi}=\xi P'+\eta R'$ and ${\tilde \eta}=\xi Q'+\eta
S'.$

\vskip 0.2cm If we identify $H_{\BR}^{(n,m)}$ with the subgroup
$$\left\{ (\xi,{\overline \xi};i\kappa)\,|\ \xi\in\BC^{(m,n)},\
\ka\in\BR^{(m,m)}\,\right\}$$ of $H_{\BC}^{(n,m)},$ we have the
following inclusion
$$G_*^J\subset SU(n,n)\ltimes H_{\BR}^{(n,m)}\subset SL(2n,\BC)\ltimes
H_{\BC}^{(n,m)}.$$ We define the mapping $\Theta:G^J\lrt G_*^J$ by
\begin{equation}\Theta
\left( \begin{pmatrix} A&B\\
C&D\end{pmatrix},(\la,\mu;\kappa) \right):=\left(\begin{pmatrix} P
& Q\\ {\overline Q} & {\overline P}
\end{pmatrix},\left(
{\frac 12}(\la+i\mu),\,{\frac 12}(\la-i\mu);\,-i{\kappa\over 2}
\right)\right),\end{equation} where $P$ and $Q$ are given by (4.13)
and (4.14). We can see that if $g_1,g_2\in G^J$, then
$\Theta(g_1g_2)=\Theta(g_1)\Theta(g_2).$

\vskip 0.2cm According to \cite[p.\,250]{YJH2002}, $G_*^J$ is of the
Harish-Chandra type\,(cf.\,\cite[p.\,118]{Sa}). Let
$$g_*=\left(\begin{pmatrix} P & Q\\
{\overline Q} & {\overline P}
\end{pmatrix},\left( \la, \mu;\,\kappa\right)\right)$$
be an element of $G_*^J.$ Since the Harish-Chandra decomposition
of an element $\begin{pmatrix} P & Q\\ R & S
\end{pmatrix}$ in $SU(n,n)$ is given by
\begin{equation*}
\begin{pmatrix} P & Q\\ R & S
\end{pmatrix}=\begin{pmatrix} I_n & QS^{-1}\\ 0 & I_n
\end{pmatrix} \begin{pmatrix} P-QS^{-1}R & 0\\ 0 & S
\end{pmatrix} \begin{pmatrix} I_n & 0\\ S^{-1}R & I_n
\end{pmatrix},
\end{equation*}
the $P_*^+$-component of the following element
$$g_*\cdot\left( \begin{pmatrix} I_n & W\\ 0 & I_n
\end{pmatrix}, (0,\eta;0)\right),\quad W\in \BD_n$$
of $SL(2n,\BC)\ltimes H_{\BC}^{(n,m)}$ is given by
\begin{equation}
\left( \begin{pmatrix} I_n & (PW+Q)(\OQ W+\OP)^{-1}
\\ 0 & I_n
\end{pmatrix},\,\left(0,\,(\eta+\la W+\mu)(\OQ W+\OP)^{-1}\,;0\right)\right).
\end{equation}

\vskip 0.2cm We can identify $\Dnm$ with the subset
\begin{equation*}
\left\{ \left( \begin{pmatrix} I_n & W\\ 0 & I_n
\end{pmatrix}, (0,\eta;0)\right)\,\Big|\ W\in\BD_n,\
\eta\in\BC^{(m,n)}\,\right\}\end{equation*} of the
complexification of $G_*^J.$ Indeed, $\Dnm$ is embedded into
$P_*^+$ given by
\begin{equation*}
P_*^+=\left\{\,\left( \begin{pmatrix} I_n & W\\ 0 & I_n
\end{pmatrix}, (0,\eta;0)\right)\,\Big|\ W=\,{}^tW\in \BC^{(n,n)},\
\eta\in\BC^{(m,n)}\ \right\}.
\end{equation*}
This is a generalization of the Harish-Chandra
embedding\,(cf.\,\cite[p.\,119]{Sa}). Then we get the {\it natural
transitive action} of $G_*^J$ on $\Dnm$ defined by
\begin{eqnarray}
& &\left(\begin{pmatrix} P & Q\\
{\overline Q} & {\overline P}
\end{pmatrix},\left( \xi, {\overline\xi};\,i\kappa\right)\right)\cdot
(W,\eta)\\
&=&\Big((PW+Q)(\OQ W+\OP)^{-1},(\eta+\xi
W+{\overline\xi})(\OQ W+\OP)^{-1}\Big),\nonumber
\end{eqnarray}

\noindent where $\begin{pmatrix} P & Q\\
{\overline Q} & {\overline P}
\end{pmatrix}\in G_*,\ \xi\in \BC^{(m,n)},\ \k\in\BR^{(m,m)}$ and
$(W,\eta)\in\Dnm.$

\vskip 0.2cm The author \cite{YJH8} proved that the action (1.2) of $G^J$
on $\Hnm$ is compatible with the action (4.17) of $G_*^J$ on $\Dnm$
through a {\it partial Cayley transform} $\Phi:\BD_{n,m}\lrt
\BH_{n,m}$ defined by
\begin{equation}
\Phi(W,\eta):=\Big(
i(I_n+W)(I_n-W)^{-1},\,2\,i\,\eta\,(I_n-W)^{-1}\Big).
\end{equation}
In other words, if $g_0\in G^J$ and $(W,\eta)\in\BD_{n,m}$,
\begin{equation}
g_0\cdot\Phi(W,\eta)=\Phi(g_*\cdot (W,\eta)),
\end{equation}
where $g_*=T_*^{-1}g_0 T_*$. $\Phi$ is a biholomorphic mapping of
$\Dnm$ onto $\Hnm$ which gives the partially bounded realization
of $\Hnm$ by $\Dnm$. The inverse of $\Phi$ is
\begin{equation*}
\Phi^{-1}(\Omega,Z)=\Big(
(\Omega-iI_n)(\Omega+iI_n)^{-1},\,Z(\Omega+iI_n)^{-1}\Big).
\end{equation*}

\vskip 0.1cm For $(W,\eta)\in \Dnm,$ we write
\begin{equation*}
(\Omega,Z):=\Phi(W,\eta).
\end{equation*}
Thus
\begin{equation}
\Om=i(I_n+W)(I_n-W)^{-1},\qquad Z=2\,i\,\eta\,(I_n-W)^{-1}.
\end{equation}
\noindent
Since
\begin{equation*}
d(I_n-W)^{-1}=(I_n-W)^{-1}dW\,(I_n-W)^{-1}
\end{equation*}
and
\begin{equation*}
I_n+(I_n+W)(I_n-W)^{-1}=2\,(I_n-W)^{-1},
\end{equation*}
we get the following formulas from (4.20)

\begin{eqnarray}
Y&=&{1 \over {2\,i}}\,(\Om-{\overline
\Om}\,)=(I_n-W)^{-1}(I_n-W\OW\,)(I_n-\OW\,)^{-1},\\
V&=&{1 \over {2\,i}}\,(Z-\OZ\,)=\eta\, (I_n-W)^{-1}+\ot\,
(I_n-\OW\,)^{-1},\\
d\Om &=& 2\,i\,(I_n-W)^{-1}dW\,(I_n-W)^{-1},\\
dZ&=&2\,i\,\Big\{ d\eta+\eta\,(I_n-W)^{-1}dW\,\Big\}(I_n-W)^{-1}.
\end{eqnarray}

Using Formulas (4.18), (4.20)-(4.24), the author \cite{YJH9} proved that
for any two positive real numbers $A$ and $B$,
the following metric $d{\tilde s}^2_{n,m;A,B}$ defined by
\begin{eqnarray*}
d{s}^2_{\Dnm;A,B}&=&4\,A\, \textrm{tr}  \Big( (I_n-W\OW)^{-1}dW(I_n-\OW W)^{-1}\bw\,\Big) \hskip 1cm\\
& &\,+\,4\,B\,\bigg\{ \textrm{tr}   \Big(
(I_n-W\OW)^{-1}\,{}^t(d\eta)\,\be\,\Big)\\
& & \quad\quad\quad
\,+\,\textrm{tr}   \Big(  (\eta\OW-{\overline\eta})(I_n-W\OW)^{-1}dW(I_n-\OW W)^{-1}\,{}^t(d\ot)\Big)\\
& & \quad\quad\quad  +\,\textrm{tr}   \Big( (\ot W-\eta)(I_n-\OW
W)^{-1}d\OW(I_n-W\OW)^{-1}\,{}^t(d\eta)\,\Big)    \\
& &\quad\quad\quad -\, \textrm{tr}   \Big( (I_n-W\OW)^{-1}\,{}^t\eta\,\eta\,
(I_n-\OW W)^{-1}\OW
dW (I_n-\OW W)^{-1}d\OW \, \Big)\\
& &\quad\quad\quad -\, \textrm{tr}   \Big( W(I_n-\OW W)^{-1}\,{}^t\ot\,\ot\,
(I_n-W\OW )^{-1}
dW (I_n-\OW W)^{-1}d\OW \,\Big)\\
& &\quad\quad\quad +\,\textrm{tr}    \Big( (I_n-W\OW)^{-1}{}^t\eta\,\ot \,(I_n-W\OW)^{-1} dW (I_n-\OW W)^{-1} d\OW\,\Big)\\
& &\quad\quad\quad +\,\textrm{tr}
\Big( (I_n-\OW)^{-1}\,{}^t\ot\,\eta\,\OW\,(I_n-W\OW)^{-1} dW (I_n-\OW W)^{-1} d\OW\,\Big)\\
& &\quad\quad\quad +\,\textrm{tr}\!  \Big( (I_n-\OW)^{-1}(I_n-W)(I_n-\OW
W)^{-1}\,{}^t\ot\,\eta\,(I_n-\OW W)^{-1}\\
& &\qquad\qquad\quad\quad \times\, (I_n-\OW)(I_n-W)^{-1}dW
(I_n-\OW W)^{-1}d\OW\,\Big)\\
& &\quad\quad\quad -\,\textrm{tr}
\Big( (I_n-W\OW)^{-1}(I_n-W)(I_n-\OW)^{-1}\,{}^t\ot\,\eta\,(I_n-W)^{-1}\\
& & \qquad\qquad\quad\quad \times\,dW (I_n-\OW
W)^{-1}d\OW\,\Big)\bigg\}
\end{eqnarray*}
\noindent is a Riemannian metric on $\Dnm$ which is invariant
under the action (4.17) of the Jacobi group $G^J_*$.

\vskip 0.2cm
We note that if $n=m=1$ and $A=B=1,$ we get
\begin{eqnarray*}
{\frac 14}\,d{ s}_{{\mathbb D}_{1,1};1,1}^2&=& { {dW\,d\OW}\over
{(1-|W|^2)^2}}\,+\,{ 1 \over {(1-|W|^2)} }\,d\eta\,d\ot\\
& &+{ {(1+|W|^2)|\eta|^2-\OW \eta^2-W\ot^2}\over {(1-|W|^2)^3} }\,dW\,d\OW\\
& & + { {\eta\OW -\ot}\over {(1-|W|^2)^2} }\,dWd\ot\,+\,{ {\ot W
-\eta}\over {(1-|W|^2)^2} }\,d\OW d\eta.
\end{eqnarray*}

\vskip 0.1cm From the formulas (4.20),\,(4.23) and (4.24), we get

\begin{equation}
\PO=\,{1\over {2\,i}}\,(I_n-W)\left[ \,{}^{{}^{{}^{{}^\text{\scriptsize $t$}}}}\!\!\!\left\{
(I_n-W)\,\PW\right\} -{}^{{}^{{}^{{}^\text{\scriptsize $t$}}}}\!\!\!\left\{ \,{}^t\eta\,{}^{{}^{{}^{{}^\text{\scriptsize $t$}}}}\!\!\!\left(
\PE\right)\right\}\,\right]
\end{equation}
and
\begin{equation}
\PZ=\,{1\over {2\,i}}\,(I_n-W)\PE.
\end{equation}


\vskip 0.3cm
Using Formulas (4.20)-(4.22), (4.25), (4.26) and Lemma 4.1, the author \cite{YJH9} proved that the following
differential operators ${\mathbb S}_1$ and ${\mathbb S}_2$ on $\Dnm$ defined by
\begin{equation*}
{\mathbb S}_1=\,\s\left( (I_n-\OW W)\PE {}^{{}^{{}^{{}^\text{\scriptsize $t$}}}}\!\!\!\left(\PEB\right)\right)
\end{equation*}
\noindent and
\begin{eqnarray*}
{\mathbb S}_2&=& \, \textrm{tr} \left( (I_n-W\OW)\,{}^{{}^{{}^{{}^\text{\scriptsize $t$}}}}\!\!\!\left(
(I_n-W\OW)\PWB\right)\PW\right)\,\\
& &   +\,\textrm{tr} \left(\,{}^t(\eta-\ot\,W)\,{}^{{}^{{}^{{}^\text{\scriptsize $t$}}}}\!\!\!\left( \PEB\right)
(I_n-\OW W)\PW  \right)\,\\
& &  +\, \textrm{tr} \left( (\ot-\eta\,\OW)\,{}^{{}^{{}^{{}^\text{\scriptsize $t$}}}}\!\!\!\left(
(I_n-W\OW)\PWB\right)\PE\right)\\
& & -\, \textrm{tr}  \left( \eta \OW
(I_n-W\OW)^{-1}\,{}^t\eta\,{}^{{}^{{}^{{}^\text{\scriptsize $t$}}}}\!\!\!\left(\PEB\right)(I_n-\OW
W)\PE\right)\\
& & -\, \textrm{tr}  \left( \ot W (I_n-\OW W)^{-1}
\,{}^t\ot\,{}^{{}^{{}^{{}^\text{\scriptsize $t$}}}}\!\!\!\left(\PEB\right)(I_n-\OW
W)\PE\right)\\
& & +\, \textrm{tr}  \left( \ot (I_n-W\OW)^{-1}{}^t\eta\,{}^{{}^{{}^{{}^\text{\scriptsize $t$}}}}\!\!\!\left(
\PEB\right)
(I_n-\OW W)\PE \right)\\
& &  +\, \textrm{tr}  \left( \eta\,\OW W (I_n-\OW
W)^{-1}\,{}^t\ot\,{}^{{}^{{}^{{}^\text{\scriptsize $t$}}}}\!\!\!\left( \PEB\right) (I_n-\OW W)\PE
\right)
\end{eqnarray*}

\noindent are invariant under the action (4.17) of $G_*^J.$ The author also proved that
\begin{equation}
\Delta_{{\mathbb D}_{n,m};A,B}:=\,{\frac 1A}\,{\mathbb S}_2\,+\,{\frac 1B}\,{\mathbb S}_1
\end{equation}
is the Laplacian of the invariant metric $ds^2_{{\mathbb D}_{n,m};A,B}$ on $\Dnm$\,(cf.\,\cite{YJH9}).

\vskip 0.5cm
\begin{proposition} The following differential operator on $\Dnm$ defined by
\begin{equation}
{\mathbb K}_{\mathbb D}=\,\det(I_n-{\overline W}W)\,\det\left( \PE {}^{{}^{{}^{{}^\text{\scriptsize $t$}}}}\!\!\!\left(
\PEB\right)\right)
\end{equation}
\noindent
is invariant under the action (4.17) of $G^J_*$ on $\Dnm$.
\end{proposition}
\vskip 0.2cm\noindent
{\it Proof.} It follows from Proposition 4.1, Formulas (4.21),\,(4.26) and the fact that the action (1.2) of $G^J$ on $\Hnm$ is
compatible with the action $(4.17)$ of $G_*^J$ on $\Dnm$ via the partial Cayley transform.
\hfill $\square$

\vskip 0.5cm
\begin{proposition} The following matrix-valued differential operator on $\Dnm$ defined by
\begin{equation}
{\mathbb T}^{\mathbb D}:=\,
{}^{{}^{{}^{{}^\text{\scriptsize $t$}}}}\!\!\!
\left( \PEB\right) (I_n-{\overline W}W) \PE
\end{equation}
\noindent
is invariant under the action (4.17) of $G^J_*$ on $\Dnm$.
\end{proposition}
\vskip 0.2cm\noindent
{\it Proof.} It follows from Proposition 4.2, Formulas (4.21),\,(4.26) and the fact that the action (1.2) of $G^J$ on $\Hnm$ is
compatible with the action $(4.17)$ of $G_*^J$ on $\Dnm$ via the partial Cayley transform.
\hfill $\square$

\vskip 0.3cm
\begin{corollary} Each $(k,l)$-entry ${\mathbb T}_{kl}^{\mathbb D}$ of ${\mathbb
T}^{\mathbb D}$ given by
\begin{equation}
{\mathbb T}_{kl}^{\mathbb D}=\sum_{i,j=1}^n
\,\left( \delta_{ij}-\sum_{r=1}^n {\overline w}_{ir}\,w_{jr}\right)
\,{{\partial^2\ \ \ \ }\over{\partial {\overline \eta}_{ki}\partial
\eta_{lj}} },\quad 1\leq k,l\leq m
\end{equation}
is a $G_*^J$-invariant differential operator on $\Dnm$.
\end{corollary}
\vskip 0.2cm\noindent
{\it Proof.} It follows immediately from Proposition 4.4.
\hfill $\square$

\vskip 0.3cm For two differential operators $D_1$ and $D_2$ on $\Hnm$ or $\Dnm$, we write
$$[D_1,D_2]:=D_1D_2-D_2D_1.$$
Then
\begin{equation}
{\mathbb M}_3=\,[{\mathbb M}_1,{\mathbb M}_2]=\,{\mathbb M}_1{\mathbb M}_2-{\mathbb M}_2{\mathbb M}_1
\end{equation}
\noindent
is an invariant differential operator of degree three on $\Hnm$ and
\begin{equation}
{\mathbb P}_{kl}=\,[{\mathbb K},{\mathbb T}_{kl}]=\,{\mathbb K}{\mathbb T}_{kl}-{\mathbb T}_{kl}{\mathbb K},
\quad 1\leq k,l\leq m
\end{equation}
\noindent
is an invariant differential operator of degree $2n+1$ on $\Hnm$.
\vskip 0.3cm
Similarly
\begin{equation}
{\mathbb S}_3=\,[{\mathbb S}_1,{\mathbb S}_2]=\,{\mathbb S}_1{\mathbb S}_2-{\mathbb S}_2{\mathbb S}_1
\end{equation}
\noindent
is an invariant differential operator of degree three on $\Dnm$ and
\begin{equation}
{\mathbb Q}_{kl}=\,[{\mathbb K}_{\mathbb D},{\mathbb T}_{kl}^{\mathbb D}]=
\,{\mathbb K}_{\mathbb D}{\mathbb T}_{kl}^{\mathbb D}-{\mathbb T}_{kl}^{\mathbb D}{\mathbb K}_{\mathbb D},
\quad 1\leq k,l\leq m
\end{equation}
\noindent
is an invariant differential operator of degree $2n+1$ on $\Dnm$.

\vskip 0.3cm
Indeed
it is very complicated and difficult at this moment to express the
generators of the algebra of all $G^J_{*}$-invariant differential
operators on $\Dnm$ explicitly.

\vskip 1cm
\end{section}

\begin{section}{{\large\bf The Case $n=m=1$}}
\setcounter{equation}{0}

\def\ddx{{{\partial^2}\over{\partial x^2}}}
\def\ddy{{{\partial^2}\over{\partial y^2}}}
\def\ddu{{{\partial^2}\over{\partial u^2}}}
\def\ddv{{{\partial^2}\over{\partial v^2}}}
\def\px{{{\partial}\over{\partial x}}}
\def\py{{{\partial}\over{\partial y}}}
\def\pu{{{\partial}\over{\partial u}}}
\def\pv{{{\partial}\over{\partial v}}}
\def\pxu{{{\partial^2\ }\over{\partial x\partial u}}}
\def\pyv{{{\partial^2\ }\over{\partial y\partial v}}}
\def\DSPR{{\Bbb D}(\SPR)}
\def\dx{{{\partial}\over{\partial x}}}
\def\dy{{{\partial}\over{\partial y}}}
\def\du{{{\partial}\over{\partial u}}}
\def\dv{{{\partial}\over{\partial v}}}
\vskip 0.2cm

\vskip 0.3cm We consider the
case $n=m=1.$ For a coordinate $(\om,z)$ in $T_{1,1}$, we write
$\om=x+i\,y,\ z=u+i\,v,\ x,y,u,v$ real. The author \cite{YJH5} proved that the algebra
$\textrm{Pol}_{1,1}^{U(1)}$ is generated by
\begin{eqnarray*}
&&q(\om,z)=\,{\frac 14}\,\omega\,{\overline \om}=\,{\frac 14}\big(x^2+y^2\big),\hskip 7cm\\
&& \xi(\om,z)=\,z\,{\overline
z}=u^2+v^2,\\
&& \phi(\om,z)=
\,{\frac 12}\,\text{Re}\,\big(z^2{\overline \om}\big)=\,{\frac 12}\,\big(u^2-v^2\big)x+uvy,\\
&& \psi(\om,z)= \,{\frac 12}\,\text{Im}\,(z^2{\overline
\om})=\,{\frac 12}\,\big( v^2-u^2\big)y+uvx.
\end{eqnarray*}

\noindent In \cite{YJH5}, using Formula (3.11) the author calculated explicitly the images
\begin{equation*}
D_1=\Theta_{1,1} (q),\quad D_2=\Theta_{1,1} (\xi),\quad D_3=\Theta_{1,1} (\phi)\quad
\textrm{and}\quad D_4=\Theta_{1,1} (\psi)
\end{equation*}
of $q,\,\xi,\,\phi$ and $\psi$ under the Halgason map $\Theta_{1,1}$.
We can show that the algebra $\BD({\mathbb
H}_{1,1})$ is generated by the following differential operators
\begin{align*} D_1=&y^2\,\left(\,{{\partial^2}\over {\partial x^2}}+
{{\partial^2}\over {\partial y^2}}\,\right)
+v^2\,\left(\,\ddu\,+\,\ddv\,\right) \\
&\ \ +2\,y\,v\,\left(\,\pxu\,+\,\pyv\,\right),
\end{align*}
$$D_2=y\left(\,{{\partial^2}\over {\partial u^2}}+
{{\partial^2}\over {\partial v^2}}\,\right),\hskip 3.74cm$$

\begin{align*}D_3=&\,y^2\,{{\partial}\over{\partial y}}
\left(\,{{\partial^2}\over{\partial u^2}}-
{{\partial^2}\over{\partial
v^2}}\,\right)\,-\,2y^2\,{{\partial^3\ \ \ \ }\over {\partial x\partial u
\partial v}}\hskip 1cm\\ &\ \ \
-\left(\, v\,{{\partial}\over{\partial v}}\,+\,1\,\right)D_2
\end{align*} and
\begin{align*} D_4=&\,y^2\,{{\partial}\over{\partial x}}\left(\,
{{\partial^2}\over{\partial v^2}}\,-\,{{\partial^2}\over {\partial
u^2}}\,\right)\,-\,2\,y^2\,{{\partial^3\ \ \ \ }\over{\partial y\partial u
\partial v}}\\ &\ \ \ \ -\,v\,{{\partial}\over{\partial u}}D_2,
\end{align*} where $\tau=x+iy$ and $z=u+iv$ with real variables
$x,y,u,v.$ Moreover, we have
\begin{align*} D_1D_2-&D_2 D_1\,=\,
2\,y^2\,\dy\left(\,\ddu\,-\,\ddv\,\right)\\ & -
4\,y^2\,{{\partial^3\ \ \ \ }\over{\partial x\partial u\partial
v}}-2\,\left(\,v\,\dv D_2+D_2\,\right).
\end{align*}\par
\noindent In particular, the algebra $\BD({\mathbb H}_{1,1})$ is
not commutative. We refer to \cite{BS, YJH5} for more detail.

\newcommand\PWE{ \frac{\partial}{\partial \overline W}}

\vskip 0.3cm
Recently Hiroyuki Ochiai \cite{Och} proved the following results.
\begin{theorem} We have the following relation
\begin{equation}
\phi^2+\psi^2 = q\,\xi^2.
\end{equation}
This relation exhausts all the relations among the generators
$q,\,\xi,\,\phi$ and $\psi$ of $\textrm{Pol}_{1,1}^{U(1)}$.
\end{theorem}

\begin{theorem} We have the following relations
\vskip 0.2cm
$(a)\ \ [D_1,D_2]=2D_3$
\vskip 0.2cm
$(b)\ \ [D_1,D_3]=2D_1D_2-2D_3$
\vskip 0.2cm
$(c)\ \ [D_2,D_3]=-D_2^2$
\vskip 0.2cm
$(d)\ \ [D_4,D_1]=0$
\vskip 0.2cm
$(e)\ \ [D_4,D_2]=0$
\vskip 0.2cm
$(f)\ \ [D_4,D_3]=0$
\vskip 0.2cm
$(g)\ \ D_3^2+D_4^2=D_2D_1D_2$
\vskip 0.3cm\noindent
These seven relations exhaust all the relations among the generators $D_1,\,D_2,\,D_3$ and $D_4$ of
$\BD({\mathbb H}_{1,1})$.
\end{theorem}

We can prove the following
\begin{theorem}
The action of $U(1)$ on $\textrm{Pol}_{1,1}^{U(1)}$ is {\it not} multiplicity-free.
\end{theorem}

\vskip 0.5cm Finally we see that for the case $n=m=1$, the seven problems proposed in Section 3 are completely solved.

\begin{remark} According to Theorem 5.2, we see that $D_4$ is a generator of the center of $\BD({\mathbb H}_{1,1})$.
We observe that the Lapalcian
$$ \Delta_{1,1;A,B}=\, {\frac 4A}\,D_1\,+\,{\frac 4B}\,D_2 \qquad (\rm{see}\ (4.8))$$
of $({\mathbb H}_{1,1},ds^2_{1,1;A,B})$ does not belong to the center of
$\BD({\mathbb H}_{1,1})$.
\end{remark}

\end{section}

\vskip 0.875cm

\begin{section}{{\large\bf The Case $n=1$ and $m$ is arbitrary}}
\setcounter{equation}{0}
\vskip 0.3cm
Conley and Raum \cite{CR} found the $2m^2+m+1$ explicit generators of $\BD({\mathbb H}_{1,m})$ and the explicit one generator
of the center of $\BD({\mathbb H}_{1,m})$. They also found the generators of
the center of the universal enveloping algebra of ${\frak U}\big({\frak g}^J\big)$ of the Jacobi Lie algebra ${\frak g}^J$. The
number of generators of the center of ${\frak U}\big({\frak g}^J\big)$ is $1+ {{m(m+1)}\over 2}.$

\vskip 0.3cm
According to Theorem 3.2, $\textrm{Pol}_{1,m}^{U(1)}$ is generated by
\begin{eqnarray}
&& q(\om,z)=\,\text{tr}(\om\,{\overline
\om}),\\
&&  \alpha_{kp}(\om,z)=
\,\text{Re}\,\big( z\,^t{\overline z} \big)_{kp}=\,\text{Re}\,(z_k{\overline z}_p),
\ 1\leq k\leq p\leq m,\\
&& \beta_{lq}(\om,z)=
\,\text{Im}\,\big( z\,^t{\overline z} \big)_{lq}=\,\text{Im}\,(z_l{\overline z}_q),
\ 1\leq l< q\leq m,\\
&& f_{kp}(\om,z)= \,\text{Re}\,(z\,{\overline
\om}\,^t\!z)_{kp}=\,\text{Re}\,({\overline \om}z_k z_p), \quad   1\leq k\leq p\leq m,\\
&& g_{kp}(\om,z)= \,\text{Im}\,(z\,
{\overline\om}\,^t\!z\,)_{kp}=\,\text{Im}\,({\overline \om}z_k z_p),\quad 1\leq k\leq p\leq m,
\end{eqnarray}
\noindent where $\om\in T_1$ and $z\in \BC^m$.

\vskip 0.3cm We let
$$\om=x+iy\in\BC\quad \textrm{and}\quad z=\,^t(z_1,\cdots,z_m)\in\BC^m \ \textrm{with}\ z_k=u_k+iv_k,\ 1\leq k\leq m,$$
where $x,y,u_1,v_1,\cdots,u_m,v_m$ are real. The invariants $q,\,\alpha_{kp},\,\beta_{lq},\,f_{kp}$ and $g_{kp}$ are
expressed in terms of $x,y,u_k,v_l\,(1\leq k,l\leq m)$ as follows:

\begin{eqnarray*}
q(\om,z)&=& x^2+y^2,\\
\alpha_{kp}(\om,z)&=& u_ku_p+v_kv_p,\quad 1\leq k\leq p\leq m,\\
\beta_{lq}(\om,z)&=& u_q v_l-u_lv_q,\quad 1\leq  l< q \leq m,\\
f_{kp}(\om,z) &=& x(u_ku_p-v_kv_p)+ y(u_kv_p+v_ku_p), \quad 1\leq k\leq p\leq m,\\
g_{kp}(\om,z)&=& x(u_kv_p+v_ku_p)-y(u_ku_p-v_kv_p), \quad 1\leq k\leq p\leq m.
\end{eqnarray*}

\vskip 0.3cm
\begin{theorem} The $1+ {{m(m+1)}\over 2}$ relations
\begin {equation}
f_{kp}^2+g_{kp}^2=q\,\alpha_{kk}\,\alpha_{pp},\quad 1\leq k\leq p\leq m
\end{equation}
exhaust all the relations among a set of generators $q,\,\alpha_{kp},\,\beta_{lq},\,f_{kp}$ and $g_{kp}$ with
$1\leq k\leq p\leq m$ and $1\leq l< q\leq m$.
\end{theorem}

\vskip 0.3cm
\begin{theorem} The action of $U(1)$ on $\textrm{Pol}_{1,m}$ is not multiplicity-free. In fact, if
$$\textrm{Pol}_{1,m}=\sum_{\sigma\in {\widehat{U(1)}} } m_\sigma \,\sigma,$$
then $m_\sigma=\infty.$
\end{theorem}

\vskip 0.3cm Problem 1, Problem 2, Problem 4, Problem 5 and Problem 7 were solved. Problem 3 can be handled.
Finally Problem 6 is unsolved in the case that $n=1$ and $m$ is arbitrary.

\end{section}

\vskip 0.875cm

\begin{section}{{\large\bf Final Remarks}}
\setcounter{equation}{0}

\newcommand\ddx{{{\partial^2}\over{\partial x^2}}}
\newcommand\ddy{{{\partial^2}\over{\partial y^2}}}
\newcommand\ddu{{{\partial^2}\over{\partial u^2}}}
\newcommand\ddv{{{\partial^2}\over{\partial v^2}}}
\newcommand\px{{{\partial}\over{\partial x}}}
\newcommand\py{{{\partial}\over{\partial y}}}
\newcommand\pu{{{\partial}\over{\partial u}}}
\newcommand\pv{{{\partial}\over{\partial v}}}
\newcommand\pxu{{{\partial^2}\over{\partial x\partial u}}}
\newcommand\pyv{{{\partial^2}\over{\partial y\partial v}}}

\vskip 0.3cm Using $G^J$-invariant differential operators on the Siegel-Jacobi
space, we introduce a noton of Maass-Jacobi forms.

\vskip 0.2cm
\begin{definition} Let
$$\Gamma_{n,m}:=Sp(n,{\mathbb Z})\ltimes H_{\mathbb Z}^{(n,m)}$$
be the discrete subgroup of $G^J$, where
$$H_{\BZ}^{(n,m)}=\left\{ (\lambda,\mu;\kappa)\in
H_{\BR}^{(n,m)}\,|\ \lambda,\mu,\kappa \ \textrm{are integral}\
\right\}.$$ A smooth function $f:\Hnm\lrt \BC$ is called a
$\textsf{Maass}$-$\textsf{Jacobi form}$ on $\Hnm$ if $f$ satisfies
the following conditions (MJ1)-(MJ3)\,:\vskip 0.1cm (MJ1)\ \ \ $f$
is invariant under $\G_{n,m}.$\par (MJ2)\ \ \ $f$ is an
eigenfunction of the Laplacian
 $\Delta_{n,m;A,B}$ (cf. Formula (4.8)).\par (MJ3)\ \ \ $f$
has a polynomial growth, that is, there exist a constant $C>0$ and
a
\par \ \ \ \ \ \ \ \ \ \ \
positive integer $N$ such that
\begin{equation*}
|f(X+iY,Z)|\leq C\,|p(Y)|^N\quad \textrm{as}\ \det
Y\longrightarrow \infty,
\end{equation*}

\ \ \ \ \ \ \ \ \ \ \ where $p(Y)$ is a polynomial in
$Y=(y_{ij}).$
\end{definition}

\begin{remark}
Let $\mathbb{D}_*$ be a commutative subalgebra of $\mathbb{D}(\Hnm)$ containing the Laplacian
$\Delta_{n,m;A,B}$.
We say that a smooth function $f:\Hnm\lrt \BC$ is a Maass-Jacobi form with respect to $\mathbb{D}_*$
if $f$ satisfies the conditions $(MJ1),\ (MJ2)_*$ and $(MJ3)$\,: the condition $(MJ2)_*$ is given by
\vskip 0.3cm\noindent $(MJ2)_*$\  $f$ is an eigenfunction of any invariant differential
operator in $\BD_*$.
\end{remark}

\vskip 0.3cm It is natural to propose the following problems.

\vskip 0.3cm\noindent {\bf {Problem\ A}\,:} Find all the
eigenfunctions of $\Delta_{n,m;A,B}.$

\vskip 0.3cm\noindent {\bf {Problem\ B}\,:} Construct Maass-Jacobi
forms.

\vskip 0.5cm If we find a {\it nice} eigenfunction $\phi$ of the Laplacian $\Delta_{n,m;A,B}$, we can construct
a Maass-Jacobi form $f_\phi$ on $\Hnm$ in the usual way defined by
\begin{equation}
f_\phi(\Omega,Z):=\,\sum_{\gamma\in \Gamma_{n,m}^\infty\backslash \Gamma_{n,m}} \phi\big( \gamma\cdot (\Omega,Z)\big),
\end{equation}
where
\begin{equation*}
\Gamma_{n,m}^\infty=\left\{ \left( \begin{pmatrix} A&B\\
C&D\end{pmatrix},(\lambda,\mu;\kappa)\right)\in \Gamma_{n,m}\,\Big|\ C=0\,\right\}
\end{equation*}
\noindent is a subgroup of $\Gamma_{n,m}.$

\vskip 0.3cm  We consider the simple case $n=m=1$ and $A=B=1$. A metric
$ds_{1,1;1,1}^2$ on $\BH_{1,1}$ given by
\begin{align*} ds^2_{1,1;1,1}\,=\,&{{y\,+\,v^2}\over
{y^3}}\,(\,dx^2\,+\,dy^2\,)\,+\, {\frac 1y}\,(\,du^2\,+\,dv^2\,)\\
&\ \ -\,{{2v}\over {y^2}}\, (\,dx\,du\,+\,dy\,dv\,)\end{align*} is
a $G^J$-invariant K{\"a}hler metric on $\BH_{1,1}$.
Its Laplacian $\Delta_{1,1;1,1}$ is given by
\begin{align*} \Delta_{1,1;1,1}\,=\,& \,y^2\,\left(\,\ddx\,+\,\ddy\,\right)\,\\
&+\, (\,y\,+\,v^2\,)\,\left(\,\ddu\,+\,\ddv\,\right)\\ &\ \
+\,2\,y\,v\,\left(\,\pxu\,+\,\pyv\,\right).\end{align*}

\vskip 0.2cm We provide some examples of eigenfunctions of
$\Delta_{1,1;1,1}$. \vskip 0.2cm (1) $h(x,y)=y^{1\over
2}K_{s-{\frac12}}(2\pi |a|y)\,e^{2\pi iax} \ (s\in \BC,$
$a\not=0\,)$ with eigenvalue $s(s-1).$ Here
$$K_s(z):={\frac12}\int^{\infty}_0 \exp\left\{-{z\over
2}(t+t^{-1})\right\}\,t^{s-1}\,dt,$$ \indent \ \ \ where
$\mathrm{Re}\,z
> 0.$ \par (2) $y^s,\ y^s x,\ y^s u\ (s\in\BC)$ with eigenvalue
$s(s-1).$
\par
 (3) $y^s v,\ y^s uv,\ y^s xv$ with eigenvalue $s(s+1).$
\par
(4) $x,\,y,\,u,\,v,\,xv,\,uv$ with eigenvalue $0$.
\par
(5) All Maass wave forms.

\vskip 0.7cm
Let $\rho$ be a rational representation of $GL(n,\BC)$ on a finite dimensional complex vector space $V_\rho$.
Let $\mathcal M\in \BR^{(m,m)}$ be a symmetric half-integral semi-positive definite matrix of degree $m$. Let
$C^\infty(\Hnm,V_\rho)$ be the algebra of all $C^\infty$ functions on $\Hnm$ with values in $V_\rho$. We define the
$|_{\rho,\mathcal M}$-slash action of $G^J$ on $C^\infty(\Hnm,V_\rho)$ as follows: If $f\in C^\infty(\Hnm,V_\rho)$,
\begin{eqnarray}
& & f|_{\rho,\mathcal M}[(M,(\lambda,\mu;\kappa))](\Om,Z) \nonumber\\
&:=&\,e^{-2\pi i \,\mathrm{tr}(\mathcal M [Z+\lambda \Omega+\mu](C\Omega+D)^{-1}C)}\cdot
e^{2\pi i \,\mathrm{tr}(\mathcal M (\lambda \Om \,^t\!\lambda\,+\,2\lambda\,^t\!Z\,+\,\kappa\,+\,\mu\,^t\!\lambda))}\\
& & \ \times\, \rho(C\Om+D)^{-1} f(M\cdot\Om,(Z+\lambda\Om+\mu)(C\Om+D)^{-1}),\nonumber
\end{eqnarray}
where $\begin{pmatrix} A&B\\
C&D\end{pmatrix}\in Sp(n,\BR)$ and $(\lambda,\mu;\kappa)\in H_\BR^{(n,m)}$. We recall the Siegel's notation
$\alpha[\beta]=\,^t\beta\alpha \beta$ for suitable matrices $\alpha$ and $\beta$. We define $\BD_{\rho,\mathcal M}$ to be
the algebra of all differential operators $D$ on $\Hnm$ satisfying the following condition
\begin{equation}
(Df)|_{\rho,\mathcal M}[g]=\,D(f|_{\rho,\mathcal M}[g])
\end{equation}
for all $f\in C^\infty(\Hnm,V_\rho)$ and for all $g\in G^J.$ We denote by ${\mathcal Z}_{\rho,\mathcal M}$ the center of
$\BD_{\rho,\mathcal M}$.

\vskip 0.3cm We define an another notion of Maass-Jacobi forms as follows.
\vskip 0.31cm
\begin{definition} A vector-valued smooth function $\phi:\Hnm\lrt V_\rho$ is called a Maass-Jacobi form on $\Hnm$ of type $\rho$
and index $\mathcal M$ if it satisfies the following conditions $(MJ1)_{\rho,\mathcal M},\ (MJ2)_{\rho,\mathcal M}$ and
$(MJ3)_{\rho,\mathcal M}$\,:
\vskip 0.1cm $(MJ1)_{\rho,\mathcal M}$\ \ \ $\phi|_{\rho,\mathcal M}[\gamma]=\phi$\ \ for all $\gamma\in\G_{n,m}.$\par
$(MJ2)_{\rho,\mathcal M}$\ \ \ $f$ is an
eigenfunction of all differential operators in the center ${\mathcal Z}_{\rho,\mathcal M}$ \par
\qquad \qquad \quad\ of $\BD_{\rho,\mathcal M}$.\par
$(MJ3)_{\rho,\mathcal M}$\ \ \ $f$
has a growth condition
$$\phi(\Om,Z)=O\Big( e^{a\det Y}\cdot e^{2\pi \textrm{tr}(\mathcal M [V]Y^{-1})} \Big)$$
\qquad \qquad \quad \ \quad as $\det Y\lrt \infty$ for some $a>0.$
\end{definition}

\vskip 0.5cm The case $n=1,\ m=1$ and $\rho=\det^k (k=0,1,2,\cdots)$ was studied by R. Bendt and R. Schmidt \cite{BS},
A. Pitale \cite{P} and K. Bringmann and O. Richter \cite{BR}. The case $n=1,\ m=$arbitrary and $\rho=\det^k (k=1,2,\cdots)$
was dealt with by C. Conley and M. Raum \cite{CR}. In \cite{CR} the authors
proved that the center ${\mathcal Z}_{\det^k,\mathcal M}$ of $\BD_{\det^k,\mathcal M}$ is the polynomial algebra with one generator
$\mathcal C^{k,\mathcal M}$, the so-called {\it Casimir} operator which is a $|_{\det^k,\mathcal M}$-slash invariant differential
operator of degree three. Bringmann and Richter \cite{BR} considered the Poincar{\'e} series ${\mathcal P}_{k,\mathcal M}$
(the case $n=m=1$)
that is a {\it harmonic} Maass-Jacobi form in the sense of Definition 7.2 and investigated its Fourier expansion and its Fourier
coefficients. Here the {\it harmonicity} of ${\mathcal P}_{k,\mathcal M}$ means that
$\mathcal C^{k,\mathcal M}{\mathcal P}_{k,\mathcal M}=0$, i.e., ${\mathcal P}_{k,\mathcal M}$ is an eigenfunction of
$\mathcal C^{k,\mathcal M}$ with zero eigenvalue. Conley and Raum \cite{CR} generalized the results in \cite{P} and \cite{BR}
to the case $n=1$ and $m$ is arbitrary.

\begin{remark}
In \cite{BCR}, Bringmann, Conley and Richter proved that the center of the algebra of differential operators invariant under the
action of the Jacobi group over a complex quadratic field is generated by two Casimir operators of degree three. They also
introduce an analogue of Kohnen's plus space for modular forms of half-integral weight over $K=\BQ (i)$, and provide a lift
from it to the space of Jacobi forms over $K$.
\end{remark}

\end{section}

\vskip 1cm

\end{document}